\documentclass[12pt]{article}

\usepackage[a4paper,includefoot,margin=2.54cm]{geometry}

\usepackage{amsmath,amsthm,amssymb}

\usepackage[utf8]{inputenc}
\usepackage[T1]{fontenc}

\usepackage{libertine,libertinust1math} 

\usepackage{mathrsfs}
\usepackage{mathtools}
\usepackage{dsfont}
\usepackage{bm}

\usepackage{graphicx}
\usepackage{float}

\usepackage[dvipsnames]{xcolor}

\usepackage[colorlinks,linkcolor=RedViolet,citecolor=PineGreen,anchorcolor=Green,urlcolor=BlueViolet]{hyperref}

\usepackage[style = authoryear, citestyle = numeric, maxnames = 20, maxcitenames = 2]{biblatex}

\DeclareFieldFormat*[article, thesis, incollection, inproceedings]{title}{#1}

\makeatletter
\input{numeric.bbx}
\makeatother

\allowdisplaybreaks

\theoremstyle{plain}  
\newtheorem{definition}{Definition}[section]

\newtheorem{lemma}[definition]{Lemma}
\newtheorem{corollary}[definition]{Corollary}

\newtheorem{proposition}[definition]{Proposition}
\theoremstyle{definition} 

\newtheorem{example}[definition]{Example}

\newtheorem{remark}[definition]{Remark}

\newcommand{\N}{\mathbb N}

\newcommand{\R}{\mathbb R}  

\newcommand{\PP}{\mathbb{P}}
\newcommand{\E}{\mathbb{E}}
\newcommand{\F}{\mathscr{F}}
\newcommand{\dt}{\mathrm{d}t}
\newcommand{\dd}{\:\mathrm{d}}
\newcommand{\e}{\mathrm{e}}
\newcommand{\oo}{\mathrm{o}}
\newcommand{\uu}{\mathrm{u}}

\usepackage{environ}
\NewEnviron{fitequation}{%
\begin{equation} \resizebox{0.91\hsize}{!}{$\displaystyle
  \BODY
$}\end{equation}
}

\NewEnviron{fitequationtwo}{%
\begin{equation} \resizebox{0.85\hsize}{!}{$\displaystyle
  \BODY
$}\end{equation}
}

\NewEnviron{fitequationtwo*}{%
\begin{equation*} \resizebox{0.94\hsize}{!}{$\displaystyle
  \BODY
$}\end{equation*}
}

\NewEnviron{fitequation*}{%
\begin{equation*} \resizebox{\hsize}{!}{$\displaystyle
  \BODY
$}\end{equation*}
}

\NewEnviron{fitalign}{
  \sbox2{\parbox{\dimexpr \wd0 + 2\wd1}%
    {\begin{align*}\BODY \end{align*}}}
  \noindent\resizebox{\textwidth}{!}{\usebox2}%
}

\numberwithin{equation}{section}
\usepackage[normalem]{ulem}

\DeclareMathOperator{\Cov}{Cov}

\usepackage{tabularray}\UseTblrLibrary{amsmath}

\begin{document}
\title{Optimal linear filtering of partially observed polynomial processes in discrete and continuous time}
\author{Jan Kallsen \quad Ivo Richert\\ Christian-Albrechts-Universität zu Kiel\footnote{Mathematisches Seminar, Christian-Albrechts-Universität zu Kiel, Heinrich-Hecht-Platz 6, 24118 Kiel, Germany, Email: kallsen@math.uni-kiel.de, richert@math.uni-kiel.de}}
\date{}
\maketitle
\begin{abstract}
This paper is devoted to filtering, smoothing, and prediction of polynomial processes that are partially observed. These problems are known to allow for an explicit solution in the simpler case of linear Gaussian state space models. The key insight underlying the present piece of research is that in linear filtering applications, polynomial processes and their discrete-time counterpart are indistinguishable from Gaussian processes sharing their first two moments. We describe the construction of these \emph{Gaussian equivalents} of polynomial processes and explicitly compute optimal linear filters, predictors and smoothers for polynomial processes in discrete and continuous time. The consideration of Gaussian equivalents also opens the door to parameter estimation and linear-quadratic optimal control in the context of polynomial processes.\\


\noindent MSC (2020) subject classification: 60G35, 62M20, 60J05, 60J25
\end{abstract}

\section{Introduction}\label{Xsec1-1}
Optimal filtering, prediction and smoothing problems constitute basic issues in applications of stochastic processes with partial observations. In mathematical statistics many parameter estimation problems require the a priori filtering of unobserved components in order to enable maximum likelihood or moment-based estimation methods. In general these problems are not easy to solve, at least if one wants to obtain optimal solutions with only moderate computational effort. For an introduction to general linear and non-linear filtering theory, see for example \textcite{bensoussan.18}, \textcite{kallianpur.13}, or \textcite{liptser.shiryaev.13}.

For linear Gaussian state space models and Gaussian Ornstein\textendash Uhlenbeck processes, however, explicit solutions to filtering problems are given by the ubiquitous K\'alm\'an filter, see for example \textcite{Anderson1979}. The situation is less obvious for non-Gaussian processes. As an example we consider the popular Heston model from mathematical finance.

\begin{example}(Discretely observed Heston model)\label{ex:heston1}
Suppose that  the logarithmic price $Y(t)$ of a stock and its stochastic squared volatility process $v(t)$ follow the Heston model
\begin{align}
\mathrm{d}v(t)&=\kappa(m-v(t))\dt+\sigma\sqrt{v(t)}\dd W^{(1)}(t),\label{e:heston1}\\
\mathrm{d}Y(t)&=\biggl(\mu+\Bigl(\delta - \frac{1}{2}\Bigr) v(t)\biggr)\dt+\sqrt{v(t)}\dd W^{(2)}(t)\label{e:heston2},
\end{align}
where $\kappa, m,\sigma \in (0, \infty)$ as well as $\mu, \delta \in \R$ denote real parameters and $W^{(1)}, W^{(2)}$ Wiener processes on some filtered probability space $(\Omega, \F, (\F_t)_{t \in \R_+}, \PP)$ with constant correlation $\rho \dt=\mathrm{d}[W^{(1)},W^{(2)}](t)$. Suppose that only $Y(t)$ is observed at discrete times $t_i\in \N$ while $v$ is unobservable. The goal is to estimate $v(t)$ from the observations $Y(t_i), i=0,1,\dots,m$ with $m \in \N$. This can be viewed as a filtering problem.

Due to their tractability, one often focuses on linear filters in practice, say of the form
\[\widehat v(t)=\alpha + \sum_{i=1}^m\gamma(i)\Delta Y(t_i)\]
for $m \in \N$ with some deterministic weights $\gamma(i)$, $i=1,\dots,m$ and $\alpha$, where we write $\Delta Y(t_i):= Y(t_i)-Y(t_{i-1})$. However, the above formulation makes little sense in the present context because $v(t)$ affects primarily the conditional variance of the returns $\Delta Y(t_i)$ rather than their mean. Instead, it seems more reasonable to consider filters of the form
\[\widehat v(t)=\alpha + \sum_{i=1}^m\gamma(i)(\Delta Y(t_i))^2\]
or possibly
\[\widehat v(t)=\alpha + \sum_{i=1}^m\beta(i)\Delta Y(t_i) + \sum_{i=1}^m\gamma(i)(\Delta Y(t_i))^2,\]
which are quadratic in the returns $\Delta Y(t_i)$ or, equivalently, linear in $(\Delta Y(t_i), (\Delta Y(t_i))^2)$,\footnote{Of course, this does not cover \emph{all} quadratic filters for $v$ because arbitrary quadratic filters could include the cross-products $\Delta Y(s)\Delta Y(t)$ across different time points $s \neq t \in \N$.} $i=1, \dots, m$. The question arises how to choose $\alpha,\beta(i)$ and $\gamma(i)$ for $i \leq m$ optimally in order to minimise the mean squared error $\E\big((\widehat v(t)-v(t))^2\big)$. If $t_m = s < t$, one speaks of an optimal prediction problem. If $s=t$, the problem is termed a filtering problem, while the case $s > t$ is commonly referred to as an optimal smoothing problem.
\label{Xenun42-1}
\end{example}

As noted above, the general filtering, prediction and smoothing problems allow for an explicit solution in discrete-time Gaussian state space models and in continuous-time Gaussian processes of Ornstein\textendash Uhlenbeck type. In this case, the optimisers are known to be deterministic linear functions of the observations and hence they coincide trivially with the optimal \emph{linear} filter, predictor, or smoother, respectively. By definition, the linear filtering problem depends only on the first two moments of the process under consideration. This yields that the K\'alm\'an filter is the optimal linear filter for all systems that share the first two moments of a Gaussian process. The same holds accordingly for the optimal prediction and smoothing problems.

On first glance, it may not seem obvious whether a setup as in Example \ref{ex:heston1} allows for an \emph{equivalent} Gaussian process in the sense of coinciding first two moments. The main insight of this paper is that such a Gaussian equivalent exists and can be determined quite easily if the process under consideration is a (possibly time-inhomogeneous) \emph{polynomial process}. This quite flexible class of Markov processes has recently attained a lot of attention in financial mathematics and modelling. It is characterised by the fact that the extended infinitesimal generator of these processes maps polynomials to polynomials of at most the same degree. This allows one to compute integer moments by solving systems of linear ordinary differential equations (ODEs), see for example \textcite{cuchiero.11}, \textcite{cuchiero.al.12}, \textcite{filipovic.larsson.16}, and \textcite{eberlein.kallsen.19}. Most importantly, polynomial processes contain the class of affine processes as a subclass, which play a predominant role in applications because of their flexibility and tractability, see \textcite{duffie.al.03}. In particular, the above Heston model fits into the framework of affine processes.

This paper is organised as follows. In Section \ref{s:pssm} we touch upon polynomial processes in continuous time and introduce the notion of a \emph{polynomial state space model} which naturally arises if a polynomial process is sampled discretely in time. These models allow for a \emph{Gaussian equivalent}, i.e.\ Gaussian processes with the same first and second moments, which are introduced in Section \ref{sec:Gauss}. This opens the door to study problems that depend only on the first two moments of the model under consideration. This applies in particular to the linear filtering, smoothing, and prediction problems, which are studied in Section \ref{s:filter}. The related issue of parameter estimation is treated in the companion paper \textcite{kallsen.richert.24b}. Further applications such as linear-quadratic control of polynomial processes are left to future research.

We generally use the notation of \textcite{eberlein.kallsen.19}. In particular, $\N = \{0, 1, \dots\}$ and $\N^* = \N \setminus \{0\}$ as well as $x^\lambda:= \prod_{j=1}^dx_j^{\lambda_j}$ for $x\in\R^d$ and $\lambda\in\N^d$. $\N^d_n$ denotes the set of multiindices $\lambda$ with $\smash{|\lambda|:= \sum_{j=1}^d\lambda_j\leq n}$. For $\lambda, \mu \in \N^d$ we write $\lambda \leq \mu$ if $\mu - \lambda \in \N^d$ and we set $\binom{\lambda}{\mu} := \prod_{j = 1}^d \binom{\lambda_j}{\mu_j}$. The multiindex $\bm{1}_j$ is defined by $(\bm{1}_j)_k=1_{\{j=k\}}$. Finally, $\bm{1}$ without a subscript denotes a multiindex of ones of suitable dimension. We write $\mathrm{I}_d$ for the identity matrix in $\R^{d \times d}$. We use the notation $k:\ell$ for  the range $k,k+1,\dots,\ell$, for example $X_{1:d}:= (X_1,\dots,X_d)$.

\section{Polynomial processes and state space models}\label{Xsec2-2}\label{s:pssm}
We first recall the properties of (time-inhomogeneous) continuous-time polynomial processes.
\subsection{Polynomial processes}\label{Xsec3-2.1}\label{su:pp}
We fix a probability space $(\Omega,\F,(\F_t)_{t\in\R_+}, \PP)$ with a right-continuous filtration $(\F_t)_{t\in\R_+}$ and an $\R^d$-valued semimartingale $X=(X(t))_{t\in\R_+}$. Following \textcite[Section VII.1]{Revuz1999}, the domain $D_{\widetilde G}$ of the extended generator of $X$ consists of those Borel-measurable functions $f: \R^d \to \R$ for which there exists a Borel-measurable function\footnote{or, strictly speaking, an equivalence class of functions $g$ modulo sets $N$ with $\int_0^\infty \PP((t, X(t)) \in N) \dd t = 0$.} $g: \R_+ \times \R^d \to \R$ such that
\[M(t):= f(X(t)) - f(X(0)) -\int_0^tg(s, X(s))\dd s, \quad t\geq0\]
is a martingale. We write $\widetilde G_t f(x) = g(t, x)$ and call $\widetilde G = (\widetilde G_t)_{t \in \R_+}$ the extended generator of $X$.

\begin{definition}\label{d:pp}
Let $\mathscr P_k(\R^d)$ denote the set of polynomials $\R^d\to\R$ of degree at most $k$. Then $X$ is called a (time-inhomogeneous) \emph{polynomial process of order $n\in\N$} if $\int_0^t \E(\lVert X(s)\rVert^n) \dd s < \infty$ for any $t \in \R_+$, if $\mathscr{P}_k(\R^d) \subset D_{\widetilde G}$ for any $k \leq n$, and if for any $\lambda \in \N^d_n$ and $f(x) = x^\lambda$ we have
\[\widetilde G_t f(x) =\sum_{\mu \in \N^d_{|\lambda|}}b^c_{\lambda,\mu}(t)x^\mu\]
for some deterministic, locally bounded coefficients $b^c_{\lambda, \mu}(t)$.  If $\smash{\widetilde G_t}$ and $\smash{B^c_n(t) := (b^c_{\lambda, \mu}(t))_{\lambda, \mu \in \N^d_n}}$ in the above definition do not depend on $t$, $X$ is called a \emph{time-homogeneous polynomial process} of order $n$.
If $X$ is a polynomial process of any order $n\in\N$, it is simply called \emph{polynomial process}. In this case we set $B^c(t) := (b^c_{\lambda, \mu}(t))_{\lambda, \mu \in \N^d}$.
\end{definition}

\begin{remark}\label{r:pp} If $n$ is even and $\E[\lVert X(0) \rVert^n] < \infty$, then the condition $\int_0^t \E(\lVert X(s) \rVert^n) \dd s < \infty$ can be omitted in the preceding definition because it follows in this case by Gr\"onwall's lemma as in (2.5) of \textcite{cuchiero.al.12}. In this case, the same argument even yields local boundedness of the map $t \mapsto \E[\lVert X(t)\rVert^n]$. Conversely, if $\int_0^t \E(\lVert X(s) \rVert^n) \dd s < \infty$ is assumed in Definition \ref{d:pp} for $t \in \R_+$, then $\E[\lVert X(0) \rVert^n] < \infty$ follows by the same Gr\"onwall argument.
\end{remark}

This family of processes has been introduced in the homogeneous setting in \textcite{cuchiero.11} and \textcite{cuchiero.al.12} and it was extended to the inhomogeneous setting in \textcite{Agoitia2020}. For properties and examples we refer to \textcite[Chapter 6]{eberlein.kallsen.19}. In particular, integer moments $\E\big(\prod_{i=1}^d X_i(t)^{k_i}\big)$ for $k_1,\dots,k_d\in\N$ can be expressed in terms of solutions to systems of linear ODEs, see \textcite[Lemma 3.3]{Agoitia2020}, or, in the homogeneous case, in terms of matrix exponentials, see \textcite[Theorem 6.26]{eberlein.kallsen.19}. We state the continuous-time moment formulas later in Eqs. \eqref{eq:poly_moments_continuous} and \eqref{e:altpssm2}.

Recall that (time-inhomogeneous) affine Markov processes in the sense of \textcite{duffie.al.03} or \textcite{filipovic.05} are polynomial of any order, given that the corresponding moments exist, see for example \textcite[Example 6.30]{eberlein.kallsen.19} in the homogeneous case. The same is true for exponential L\'evy processes or, more generally, generalised Ornstein\textendash Uhlenbeck processes in the sense of \textcite[Section 3.9]{eberlein.kallsen.19}.

\begin{remark}\label{r:pp1}
Note that $X$ is polynomial of order 1 if and only if $\smash{\int_0^t \E(\lVert X(s)\rVert) \dd s  < \infty}$ for any $t \in \R_+$ and if, for some measurable, locally bounded coefficients $a^c(t) \in \R^d$, $A^c(t) \in \R^{d \times d}$,
\begin{equation}\label{eq:pp1}
\mathrm{d}X(t)=(a^c(t) + A^c(t)X(t))\dt + \mathrm{d}M(t)
\end{equation}
holds for $t \in \R_+$, where the process $M = (M(t))_{t \in \R_+}$ is an $\R^d$-valued martingale.
\end{remark}
Polynomial processes of order $n$ can be characterised in terms of those of order 1:

\begin{lemma}\label{l:pp}
The following statements are equivalent:
\begin{enumerate}
\item $X$ is a polynomial process of order $n$.
\item For each $m\in \{1,\dots,n\}$, the $\R^{\widetilde m}$-valued processes $(X^\lambda)_{\lambda\in\N^d_m\setminus\{0\}}$ are polynomial processes of order 1, where $\widetilde m = \widetilde m(m)$ is defined by $\widetilde m:= \binom{d+m}{m}-1$.
\item For each $m\in \{1,\dots,n\}$, the $\R^{(1+d)^m}$-valued processes $(1,X)^{\otimes m}$ are polynomial processes of order 1, where $\otimes m$ denotes the $m$-fold Kronecker product of a vector.
\end{enumerate}
\end{lemma}

\begin{proof}
Any of the three statements above means that the drift rate $\widetilde G_t f(X(t))$ of monomials $f(X) = X^\lambda$ is a polynomial in $X(t)$ of at most the same degree for each multiindex $\lambda \in \N^d_n$.
\end{proof}

\begin{remark}Let $\smash{B^c(t)=(b^c_{\lambda,\mu}(t))_{\lambda,\mu\in\N^d}}$ denote the matrix belonging to $X$ in the sense of Definition \ref{d:pp}. Once a specific ordering $(\lambda_1, \dots, \lambda_{\widetilde m})$ of $\N^d_m\setminus \{0\}$ has been chosen, the coefficients $a^c(t) \in \R^{\widetilde m}$ and $A^c(t)\in \R^{\widetilde m\times \widetilde m}$ in Remark \ref{r:pp1} for  $(X^\lambda)_{\lambda\in\N^d_m\setminus\{0\}}$ can be written in the form
\begin{align*}
a^c_j(t) &= b^c_{\lambda_j, 0}(t), \quad A^c_{j, k}(t) = b^c_{\lambda_j, \lambda_k}(t).
\end{align*}
The components of the processes  $(1,X)^{\otimes m}$ are of the form $X^\lambda$ for some $\lambda$ satisfying $|\lambda|\leq m$. Since the same exponent $\lambda\in\N^d$ appears several times, the coefficients $a^c$ and $A^c$ in Remark \ref{r:pp1} are not uniquely specified for $(1,X)^{\otimes m}$. \mbox{The non-unique $a^c$ and $A^c$ may be chosen as follows:}
\begin{align*}
a^c_j(t)&\equiv 0,\\
A^c_{j,k}(t)&=
\begin{cases}
0 & \text{if } k-1\notin I\\
b^c_{\lambda(j-1),\:\lambda(k-1)}(t) & \text{if } k-1\in I,
\end{cases}
\end{align*}

\noindent where $\lambda(j):= (\lambda_1,\dots,\lambda_d)$ with $\lambda_i$ being the number of digits $i\in\{0,\dots,d\}$ in the base-$(1+d)$ representation of $j=(a_\ell\dots a_1a_0)_{1+d}$ ( $\lambda(38)=\lambda(1102_3)=(2,1)$ for $d=2$) and $I\subset\N$ denotes the set of numbers whose digits in the base-$(1+d)$ representation are increasing, i.e., for $d = 2$,
\[I=\{0_3,1_3,2_3,11_3,12_3,22_3,111_3,112_3,122_3,222_3,1111_3,1112_3,1122_3,1222_3,2222_3,\dots\}.\]
\end{remark}
Existing polynomial processes can be used to construct new ones. We consider two cases.
\begin{lemma}\label{l:vec}
Let $m\in\N^*$. If $X$ is polynomial of any order, the same is true for $(X^\lambda)_{\lambda\in\N^d_m\setminus\{0\}}$.
\end{lemma}

\begin{proof}
This follows easily from the definition.
\end{proof}

The next example contains generalised Ornstein\textendash Uhlenbeck processes in the sense of \textcite{Behme2012} as a special case. It also entails that the K\'alm\'an\textendash Bucy filter from Section \ref{su: KB} for a polynomial process of order 2 preserves the polynomial property.

\begin{lemma}\label{l:gou}
Suppose that $X$ is polynomial of order $n \geq 2$, let $Y$ be an independent $\R^{k\times k}$-valued L\'evy process with $\E(\lVert Y(1)\rVert^n) < \infty$, and let $c(t)\in\R^k$, $C(t)\in\R^{k\times d}$ be locally bounded. Moreover, let $Z$ be an $\R^k$-valued adapted process with $\int_0^t \E(\lVert Z(s)\rVert^n) \dd s  < \infty$ for $t \in \R_+$ that solves
\[ \mathrm{d}Z(t)= \mathrm{d}Y(t) Z(t-) + c(t) \dt + C(t) \dd X(t).\]
If additionally $\E\bigl(\int \lVert \xi\rVert^n K_X(t, \mathrm{d}\xi)\bigr) < \infty$, $t \in \R_+$, where $K_X$ denotes the third differential characteristic of $X$, then $(X,Z)$ is an $\R^{d+k}$-valued polynomial process of order $n$. If $c(t) \equiv c$, $C(t) \equiv C$ and $X$ is time-homogeneous, then $(X, Z)$ is a time-homogeneous polynomial process of order $n$.
\end{lemma}

\begin{proof}
For ease of notation, we consider the case $d=1=k$. The general statement follows along the same lines. Denote by $I(t):= t$ the identity process and by $\psi^I,\psi^X,\psi^Y$ the characteristic exponents of $I$, $X$, $Y$ in the sense of \textcite[Theorem 2.17 and Definition 4.3]{eberlein.kallsen.19}. Since $I$, $X$, $Y$ are independent, \mbox{the characteristic exponent of $(I, X, Y)$ is given by}
\[\psi^{(I, X, Y)}(t,(u_1,u_2,u_3))=\psi^I(u_1)+\psi^X(t, u_2)+\psi^Y(u_3)\]
for $u = (u_1, u_2, u_3) \in \R^3$ with $\psi^I(u_1)=iu_1$ for $u_1 \in \R$. Observe further that
\[\binom{X}{Z}=\binom{X(0)}{Z(0)}
+ \begin{pmatrix}
0 & 1 & 0\\
c & C & Z_-
\end{pmatrix}
\bullet
(I, X, Y)^\top,\]
where the dot denotes stochastic integration. Therefore, by \textcite[Proposition 4.8]{eberlein.kallsen.19}, the characteristic exponent of the $\R^{k+d}$-valued process $(X,Z)$ is of the form
\begin{align*}
\psi^{(X,Z)}(t,(u_1,u_2))
&=\psi^{(I,X,Y)}\left(t, \:
\begin{pmatrix}
0 & 1 & 0\\
c & C & Z(t-)
\end{pmatrix}^\top
\binom{u_1}{u_2}\right)\\
&=ic(t)u_2 + \psi^X(t,u_1+C(t)u_2) + \psi^Y(Z(t-)u_2).
\end{align*}
Let $\lambda=(\lambda_1,\lambda_2)\in\N^2_n$. By the arguments in the proof of Theorem 5.2 in \textcite{Agoitia2020}, $X$ admits differential characteristics $(b_X, c_X, K_X)$, satisfying $\int \lVert \xi \rVert^n K_X(t, \mathrm{d}\xi) < \infty$ by assumption. Therefore, $\psi^X(t, u)$ is $n$ times differentiable in $u$ and by the same arguments as in the proof of  Theorem 6.25 in \textcite{eberlein.kallsen.19}, we have that
\[\frac{\partial^{\lambda} \psi^X(t,u_1+C(t)u_2)}{\partial (u_1,u_2)^{\lambda}}\bigg|_{u=0}
= \frac{\partial^{|\lambda|}\psi^X(t,u_3)}{\partial u_3^{|\lambda|}}\bigg|_{u_3=0}C(t)^{\lambda_2}\]
is a polynomial of order at most $|\lambda|$ in $X(t-)$ because $X$ is polynomial of order $n$. The coefficients of this polynomial are locally bounded in $t$ because they are linear combinations of the coefficients $b^c_{\lambda, \mu}$ corresponding to $X$ (see equation (6.52) in \textcite{eberlein.kallsen.19} and the corresponding argument in the proof of Theorem 5.3 in \textcite{Agoitia2020}). By the integrability assumption on $Y$, $\psi^Y$ is also $n$ times differentiable and we have
\[\frac{\partial^{\lambda_2} \psi^Y(Z(t-)u_2)}{\partial u_2^{\lambda_2}}\bigg|_{u_2=0}=
\frac{\partial^{\lambda_2} \psi^Y(u_2)}{\partial u_2^{\lambda_2}}\bigg|_{u_2=0}Z(t-)^{\lambda_2},\]
which implies that it is a polynomial of order at most $\lambda_2$ in $Z(t-)$. Altogether, we have that ${\partial\psi^{(X,Z)}(t,u)/\partial u^\lambda}|_{u=0}$ is a polynomial of order at most $|\lambda|$ in $(X(t-),Z(t-))$. Once we show
\begin{equation}\label{e:stern}
\E\biggl(\int \| \xi \|^n K_{(X, Z)}(t, \mathrm{d}\xi)\biggr) < \infty,
\end{equation}
where $K_{(X, Z)}$ is the third differential characteristic of $(X, Z)$, then we can apply the proof of Theorem 5.3 of \textcite{Agoitia2020}, which yields that the extended generator of the process $(X, Z)$ at $t$ maps $x^\lambda$ to a polynomial in $x$ of degree at most $|\lambda|$ with locally bounded coefficients. \eqref{e:stern} however holds by assumption because $\int \lVert \xi \rVert^n K_{(X, Z)}(t, \mathrm{d}\xi) \leq C[\int \lVert \xi \rVert^n K_{X}(t, \mathrm{d}\xi) + \lVert Z(t-)\rVert^n]$ for some constant $C > 0$, which can be deduced by \textcite[Proposition 4.8]{eberlein.kallsen.19}. The map $t \mapsto \E(\lVert (X(t), Z(t))\rVert^n)$ is locally integrable because $X$ is polynomial of order $n$ and by assumption on $Z$. Hence $(X, Z)$ is a polynomial process of order $n$. The final statement for the time-homogeneous case follows along the same lines.
\end{proof}

\subsection{Polynomial state space models}\label{Xsec4-2.2}\label{su:pssm}
We turn now to the discrete-time counterpart of polynomial processes. To this end, we fix a filtered probability space $(\Omega,\F,(\F_t)_{t\in\N},\PP)$ and an $\R^d$-valued adapted process $X=(X(t))_{t\in\N}$.
\begin{definition}\label{d:poly1}
$X$ is a \emph{polynomial state space model of order 1} if \ $\E(\lVert X(t)\rVert)<\infty$, $t\in\N$, and
\begin{equation}\label{e:definition_time_discrete_poly_ssm}
X(t)=a(t)+A(t)X(t-1) + N(t),\quad t \in \N^*,
\end{equation}
for some deterministic $a(t)\in\R^d$, $A(t)\in\R^{d\times d}$ and some  $\R^d$-valued martingale difference sequence $N$ satisfying $\E(N(t)|\F_{t-1}) = 0$ for any $t\in\N^*$. We say that $X$ has \emph{second moments} if $\E(\lVert X(t)\rVert^2)<\infty$, $t\in\N$. \mbox{Moreover, we call $X$ \emph{time-homogeneous} if $a(t),A(t)$ do not depend on $t$.}
\end{definition}

In line with the definition of continuous-time polynomial processes in Section \ref{su:pp} and in particular with Lemma \ref{l:pp}, polynomial state space models can be defined for any order.
\begin{definition}\label{d:pssm}
We call $X$ \emph{(time-homogeneous) polynomial state space model of order $n$} if $(X^\lambda)_{\lambda\in\N^d_m\setminus\{0\}}$ or, equivalently, $(1,X)^{\otimes m}$ is a (time-homogeneous) polynomial state space model of order 1 for $m=1,\dots,n$. Moreover, we call $X$ a \emph{(time-homogeneous) polynomial state space model} if it is a (time-homogeneous) polynomial state space model of any order $n\in\N$.
\end{definition}
The following characterisation corresponds to Definition \ref{d:pp}.
\begin{lemma}\label{l:altpssm}
$X$ is a polynomial state space model of order $n$ if and only if for $t \in \N^*$
\begin{equation}\label{e:altpssm}
\E\bigl(X(t)^\lambda\big|\F_{t-1}\bigr) = \sum_{\mu\in\N^d_n} b_{\lambda,\mu}(t)X(t-1)^\mu
\end{equation}
for any multiindex $\lambda\in\N^d_n$ for some $b_{\lambda,\mu}(t)\in\R$ such that $b_{\lambda,\mu}(t)=0$ if $|\mu|>|\lambda|$. In this case, $X$ is time-homogeneous if the coefficients $b_{\lambda,\mu}(t)$ do not depend on $t$. Also, the coefficients $a(t),A(t)$ in \eqref{e:definition_time_discrete_poly_ssm} for $(X^\lambda)_{\lambda\in\N^d_m\setminus\{0\}}$ are given by $a_j(t)=b_{\lambda_j,0}(t)$ and $A_{i, j}(t)=b_{\lambda_i,\lambda_j}(t)$, $i, j \in \{1, \dots, \widetilde m\}$.
\end{lemma}

\begin{proof}
This follows directly from the definition.
\end{proof}

As in continuous time, conditional and unconditional moments can be computed explicitly.
\begin{proposition}[Moment formula]\label{p:moments}
If $X$ is a polynomial state space model of order $n$, then
\begin{equation}\label{e:moments}
\E\bigl((X^\lambda(t))_{\lambda\in\N^d_n\setminus\{0\}}\big|\F_s\bigr) =\Big(\prod_{r=s+1}^t A(r) \Big)(X^\lambda(s))_{\lambda\in\N^d_n\setminus\{0\}}+\sum_{r=s+1}^{t}\Big( \prod_{u=r+1}^t A(u)\Big)a(r)
\end{equation}
and
\begin{align}\label{e:moments2}
\E\bigl((X^\lambda(t))_{\lambda\in\N^d_n}\big|\F_s\bigr)
=&\Big(\prod_{r=s+1}^t B(r)\Big)(X^\lambda(s))_{\lambda\in\N^d_n}
\end{align}
for $t \in \N$ and $s \leq t$, where $a(t), A(t)$ are the coefficients in \eqref{e:definition_time_discrete_poly_ssm} for $(X^\lambda)_{\lambda\in\N^d_n\setminus\{0\}}$ instead of $X$ and $B(t):= (b_{\lambda,\mu}(t))_{\lambda,\mu\in\N^d_n}$ is the matrix from \eqref{e:altpssm} for $X$. A parallel statement holds for $(1,X)^{\otimes n}$.
\label{Xenun20-1}
\end{proposition}

\begin{proof}
This follows by induction on $t-s$.
\end{proof}

\begin{remark}\label{rem: time-hom}
In the time-homogeneous case where the coefficient vector and matrices $a:= a(t), A:= A(t)$ and $ B:= B(t)$ do not depend on $t$, Eqs. \eqref{e:moments}, \eqref{e:moments2} boil down to
\begin{align*}
\E\bigl((X^\lambda(t))_{\lambda\in\N^d_n\setminus\{0\}}\big|\F_s\bigr)
&=A^{t-s}(X^\lambda(s))_{\lambda\in\N^d_n\setminus\{0\}}+\sum_{r=s+1}^{t}A^{t-r}a,\\
\E\bigl((X^\lambda(t))_{\lambda\in\N^d_n}\big|\F_s\bigr)
&=B^{t-s}(X^\lambda(s))_{\lambda\in\N^d_n}.
\end{align*}

\pagebreak
\noindent If $\rho(A) < 1$, i.e.\ the spectral radius of $A$ is less than 1, these moments converge to
\begin{align}\label{e:nocite}
\lim_{t\to\infty}\E\bigl((X^\lambda(t))_{\lambda\in\N^d_n\setminus\{0\}}\big|\F_s\bigr)
&=\sum_{r=0}^\infty A^ra=(\mathrm{I}-A)^{-1}a
\end{align}
almost surely and also in $L^1$, where $\mathrm{I}$ denotes the identity matrix of suitable size.
\end{remark}

Lemmas \ref{l:vec} and \ref{l:gou} hold accordingly in the present setup:
\begin{lemma}\label{l:vec2}
Let $m\in\N^*$. If $X$ is a polynomial state space model, the same holds for the process $(X^\lambda(t))_{\lambda\in\N^d_m\setminus\{0\}}$.
\end{lemma}

\begin{proof}
This follows easily from the definition.
\end{proof}

\begin{proposition}\label{p:gou2}
Suppose that $X$ is a polynomial state space model of order $n$ and let $(Y(t))_{t \in \N^*}$ be an independent adapted sequence of $\R^{k\times k}$-valued random variables with $\E[\lVert Y(t)\rVert^n] < \infty$ such that $Y(t)$ is independent of $\F_{t-1}$. Moreover, let $c(t)\in\R^{k}$, $C(t)\in\R^{k\times d}$, $t \in \N^*$, be deterministic and $Z(0)$ an $\mathscr F_0$-measurable $\R^k$-valued random variable with $\E[\lVert Z(0)\rVert^n] < \infty$. Let the $\R^k$-valued process $Z$ be given by
\begin{equation}\label{e:mitt}
Z(t) =  Y(t)Z(t-1) + c(t) + C(t)X(t)
\end{equation}
for $t \in \N^*$, or, alternatively, given by
\begin{equation}\label{e:mitt-1}
Z(t) =  Y(t)Z(t-1) + c(t) + C(t)X(t-1)
\end{equation}
for $t \in \N^*$. Then the following two statements hold true:
\begin{enumerate}
\item The process $(X,Z)$ is an $\R^{d+k}$-valued polynomial state space model of order $n$.
\item Let $\lambda = (\lambda_1, \lambda_2), \mu = (\mu_1, \mu_2) \in \N^{d + k}_n$. The moment formula for $(X, Z)$ reads as follows.
\begin{align*}
\E\bigl((X, Z)(t)^\lambda |\F_{t-1}\bigr) = \sum_{\mu \in \N^{d+k}_{|\lambda|}}
\widebar b_{\lambda, \mu}(t) \;(X, Z)(t-1)^\mu,
\end{align*}
where the coefficients $\widebar b_{\lambda, \mu}(t)$, $t \in \N^*$, corresponding to the process $(X, Z)$ are given by
\begin{align*}
&\widebar b_{\lambda, \mu}(t) := \begin{cases}
\displaystyle \sum_{\nu \in \N^d_n} \sum_{\eta\leq\lambda_2}  \textstyle \binom{\lambda_2}{\eta}
S_{\nu,\eta}(t)\widetilde S_{\mu_2,\lambda_2-\eta}(t)
b_{\lambda_1+\nu,\mu_1}(t)
& \text{ for \eqref{e:mitt}},\\
\displaystyle\sum_{\nu\leq\mu_1} \sum_{\eta\leq\lambda_2}  \textstyle \binom{\lambda_2}{\eta}
S_{\nu,\eta}(t)\widetilde S_{\mu_2,\lambda_2-\eta}(t)
b_{\lambda_1,\mu_1-\nu}(t)& \text{ for \eqref{e:mitt-1}}.
\end{cases}
\end{align*}
Here, let $M_{\nu, \eta} := \{\alpha \in \N_n^{k \times d}: \alpha \bm{1} = \eta, \bm{1}^\top \alpha = \nu\}$, $N_{\nu, \eta} := \{\alpha \in \N_n^{k \times k}: \alpha \bm{1} \leq \eta, \bm{1}^\top \alpha = \nu\}$, and
\begin{align*}
S_{\nu, \eta}(t) &:= \sum_{\alpha \in M_{\nu, \eta}} \prod_{j=1}^k \textstyle \binom{\eta_j}{\alpha_{j, \cdot}} C(t)_{j, \cdot}^{\alpha_{j, \cdot}}, \quad \displaystyle\widetilde S_{\nu, \eta}(t) := \sum_{\alpha \in N_{\nu, \eta}} \E\Big(\prod_{j=1}^k \textstyle \binom{\eta_j}{(\alpha_{j, \cdot}, \:\eta_j - |\alpha_{j, \cdot}|)} Y(t)_{j, \cdot}^{\alpha_{j, \cdot}} c(t)_j^{\eta_j - |\alpha_{j, \cdot}|}\Big).
\end{align*}
Moreover, $(b_{\lambda,\mu}(t))_{\lambda,\mu\in\N^d_n}$ denotes the matrix of coefficients in the sense of \eqref{e:altpssm} for $X$.
\end{enumerate}
\end{proposition}

\begin{proof}

In the calculations below we start with the case \eqref{e:mitt} because the calculations for \eqref{e:mitt-1} are similar. First suppose $c(t) \equiv 0$. By the binomial theorem it follows that
\begin{fitequation*}
    X(t)^{\lambda_1} Z(t)^{\lambda_2} = X(t)^{\lambda_1} \big( Y(t) Z(t-1) + C(t) X(t)\big)^{\lambda_2} = X(t)^{\lambda_1} \sum_{\eta \leq \lambda_2} \binom{\lambda_2}{\eta} \big(C(t)X(t)\big)^{\lambda_2 - \eta} \big( Y(t) Z(t-1)\big)^{\eta}
\end{fitequation*}
for $\lambda = (\lambda_1, \lambda_2) \in \N^{d+ k}_n$.
The binomial theorem applied to the product $(Y(t) Z(t-1))^\eta$  implies
\begin{align*}
    \big( Y(t) Z(t-1)\big)^{\eta} &= \prod_{j=1}^k \sum_{\substack{\alpha \in \N^k_n \\ |\alpha| = \eta_j}} \binom{\eta_j}{\alpha} Y(t)_{j,\cdot}^\alpha Z(t-1)^\alpha\\
    &= \sum_{\substack{\alpha_1, \dots, \alpha_k \in \N^k_n \\ |\alpha_i| = \eta_i}} \Big[\prod_{j=1}^k \binom{\eta_j}{\alpha_j} Y(t)_{j, \cdot}^{\alpha_j}\Big] Z(t-1)^{\alpha_1 + \dots + \alpha_k} \\
    &= \sum_{\substack{\mu_2 \in \N^k_n \\ |\mu_2| = |\eta|}} \sum_{\alpha \in \widehat N_{\mu_2, \eta}}\Big[\prod_{j=1}^k \binom{\eta_j}{\alpha_j} Y(t)_{j, \cdot}^{\alpha_j}\Big] Z(t-1)^{\mu_2} =: \sum_{\substack{\mu_2 \in \N^k_n \\ |\mu_2| = |\eta|}} \widehat S_{\mu_2, \eta}(t) Z(t-1)^{\mu_2},
\end{align*}
where we set $\widehat N_{\mu_2, \eta} := \{\alpha \in \N^{k \times k}_n:\: \alpha\bm{1} = \eta, \bm{1}^\top \alpha = \mu_2\}$. Inserting this into the above yields
\begin{align*}
    X(t)^{\lambda_1} Z(t)^{\lambda_2} &=  X(t)^{\lambda_1} \sum_{\eta \leq \lambda_2} \sum_{\substack{\mu_2 \in \N^k_n \\ |\mu_2| = |\eta|}} \binom{\lambda_2}{\eta} \big(C(t)X(t)\big)^{\lambda_2 - \eta}  \widehat S_{\mu_2, \eta}(t) Z(t-1)^{\mu_2} \\ &= X(t)^{\lambda_1} \sum_{\mu_2 \in \N^k_{|\lambda_2|}} \underbrace{\sum_{\substack{\eta \leq \lambda_2 \\ |\eta| = |\mu_2|}}  \binom{\lambda_2}{\eta} \big(C(t)X(t)\big)^{\lambda_2 - \eta}  \widehat S_{\mu_2, \eta}(t)}_{=: \gamma_{\lambda_2, \mu_2}(t)} Z(t-1)^{\mu_2}.
\end{align*}
Analogously to the calculation for $\widehat S_{\mu_2, \eta}(t)$ above, we obtain for the  $\gamma_{\lambda_2, \mu_2}(t)$, $t \in \N$, that

\begin{align*}
    \gamma_{\lambda_2, \mu_2}(t) &= \sum_{\eta \leq \lambda_2} \binom{\lambda_2}{\eta} \widehat S_{\mu_2, \eta}(t) \bm{1}_{\{|\eta| = |\mu_2|\}} \big(C(t)X(t)\big)^{\lambda_2 - \eta}\\ &\stackrel{\lambda_2 - \eta \rightarrow \eta}{=} \sum_{\eta \leq \lambda_2} \binom{\lambda_2}{\eta} \widehat S_{\mu_2, \lambda_2 - \eta}(t) \bm{1}_{\{|\eta| = |\lambda_2| - |\mu_2|\}} \big(C(t)X(t)\big)^{\eta} \\
    &= \sum_{\nu \in \N^d_{|\lambda_2|}} \sum_{\substack{\eta \leq \lambda_2 \\ |\eta| = |\nu|}} \binom{\lambda_2}{\eta} \widehat S_{\mu_2, \lambda_2 - \eta}(t) \bm{1}_{\{|\eta| = |\lambda_2| - |\mu_2|\}} S_{\nu, \eta}(t) X(t)^\nu \\
    &= \sum_{\substack{\nu \in \N^d \\ |\nu| = |\lambda_2| - |\mu_2|}} \sum_{\substack{\eta \leq \lambda_2 \\ |\eta| = |\nu|}} \binom{\lambda_2}{\eta} \widehat S_{\mu_2, \lambda_2 - \eta}(t) S_{\nu, \eta}(t) X(t)^\nu,
\end{align*}

\noindent where $S_{\nu, \eta}(t)$ is defined as in the proposition. In the above sums, we can omit the subscripts $|\eta| = |\nu|$ and $|\nu| = |\lambda_2| - |\mu_2|$ because $S_{\nu, \eta}(t) = 0$ unless $|\eta| = |\nu|$ and $\widehat S_{\mu_2, \lambda_2 - \eta}(t) = 0$ unless $|\lambda_2| - |\mu_2| = |\eta|$ $(= |\nu|)$. Inserting this into the decomposition for $X(t)^{\lambda_1}Z(t)^{\lambda_2}$ yields
\begin{equation*}
     X(t)^{\lambda_1} Z(t)^{\lambda_2} = \sum_{\mu_2 \in \N^k_{|\lambda_2|}} \sum_{\nu \in \N^d}  \sum_{\eta \leq \lambda_2} \binom{\lambda_2}{\eta} \widehat S_{\mu_2, \lambda_2 - \eta}(t) S_{\nu, \eta}(t) Z(t-1)^{\mu_2} \begin{cases}
         X(t)^{\lambda_1 + \nu} & \text{ for } \eqref{e:mitt}, \\
         X(t)^{\lambda_1} X(t-1)^{\nu}& \text{ for } \eqref{e:mitt-1}.
     \end{cases}
\end{equation*}
Since $Y$ is independent of $X$ and $Y(t)$ is independent of $\F_{t-1}$, this implies that $(X, Z)$ is polynomial of order $n$. For the case $\eqref{e:mitt}$ we apply $\E(X(t)^{\lambda_1 + \nu} |\F_{t-1}) = \sum_{\mu \in \N^d_n} b_{\lambda_1 + \nu, \mu}(t) X(t-1)^\mu$, while for the case \eqref{e:mitt-1} we apply $\E(X(t)^{\lambda_1} X(t-1)^\nu |\F_{t-1}) = \sum_{\mu \in \N^d} b_{\lambda_1, \mu}(t) X(t-1)^{\mu + \nu} = \sum_{\mu \geq \nu} b_{\lambda_1, \mu - \nu}(t) X(t-1)^\mu$ to obtain $\E((X, Z)(t)^\lambda |\F_{t-1}) = \sum_{\mu \in \N^{d+k}_{|\lambda|}} \overline{b}_{\lambda, \mu}(t) (X, Z)(t-1)^\mu$ with
\[\overline{b}_{\lambda, \mu}(t) = \begin{cases}
    \sum_{\nu \in \N^d_n} \sum_{\eta \leq \lambda_2} \binom{\lambda_2}{\eta} \widehat S_{\mu_2, \lambda_2 - \eta}(t) S_{\nu, \eta}(t) b_{\lambda_1 + \nu, \mu_1}(t) & \text{ for } \eqref{e:mitt}, \\
    \sum_{\nu \leq \mu_1} \sum_{\eta \leq \lambda_2} \binom{\lambda_2}{\eta} \widehat S_{\mu_2, \lambda_2 - \eta}(t) S_{\nu, \eta}(t) b_{\lambda_1, \mu_1 - \nu}(t) & \text{ for } \eqref{e:mitt-1}.
\end{cases}\]
Incorporating the case $c(t) \neq 0$ leads to similar calculations by considering the process $\widetilde Z(t) := (Z(t), 1)$, which then fulfils the homogeneous equation $\widetilde Z(t) = \widetilde Y(t) \widetilde Z(t-1) + \widetilde C(t) X(t)$ in the case \eqref{e:mitt} with the obvious modifications in the case \eqref{e:mitt-1} for some matrices $\widetilde C(t) \in \R^{(k+1)\times d}$ and $\widetilde Y(t) \in \R^{(k+1)\times (k+1)}$. Hence the above calculations can be applied to $(X, \widetilde Z)$.
\end{proof}

Dynamic equations of the form \eqref{e:mitt} and \eqref{e:mitt-1} occur often in filtering applications. Indeed, the popular K\'alm\'an filter recursions stated in Proposition \ref{p:Ka-lman} are exactly of the form \eqref{e:mitt} for the optimal filter and of the form \eqref{e:mitt-1} for the optimal predictor. Therefore, Proposition \ref{p:gou2} states that filtered polynomial state space models as in Section \ref{su: Kal} stay polynomial state space models and yields an explicit moment formula for them.

Polynomial state space models occur naturally if polynomial processes are discretely sampled. The following proposition provides simple formulas to cast the coefficients of polynomial processes from Definition \ref{d:pp} and Remark \ref{r:pp1} into the corresponding coefficients for polynomial state space models from \eqref{e:definition_time_discrete_poly_ssm} and \eqref{e:altpssm} if the former are sampled at $t_0 \leq t_1 \leq t_2 \leq \dots$.
\begin{proposition}\label{l:restriction}
If $X = (X(t))_{t \in \R_+}$ is a polynomial process of order $n$, then $(X(t_k))_{k\in\N}$ is a polynomial state space model of order $n$ with respect to $(\F_{t_k})_{k \in \N}$ that is homogeneous if $X$ is homogeneous and if $\Delta t:= t_{k + 1} - t_{k}$ does not depend on $k$. Moreover:
\begin{enumerate}
\item Let $a^c(t)$ and $A^c(t)$ denote the coefficients from Remark \ref{r:pp1} for $(X^\lambda(t))_{\lambda\in\N^d_n\setminus\{0\}}$, which satisfy
\begin{equation}\label{eq:longdynamics}
\mathrm{d}(X^\lambda(t))_{\lambda\in\N^d_n\setminus\{0\}}=\bigl(a^c(t)+A^c(t) \:(X^\lambda(t))_{\lambda\in\N^d_n\setminus\{0\}}\bigr)\dt+\mathrm{d}M(t)
\end{equation}
for some martingale $M$. Let $F(t)$ denote the unique solution to $\mathrm{d}F(t) = A^c(t) F(t) \dd t$ with $F(0) = \mathrm{I}$. Set $N(t_k) := F(t_k) \int_{t_{k-1}}^{t_k} F(s)^{-1} \dd M(s)$. Then $\E(N(t_k)|\F_{t_{k-1}}) = 0$ and
\begin{equation}\label{eq:poly_moments_continuous}
(X^\lambda(t_k))_{\lambda\in\N^d_n\setminus\{0\}}= a(k) + A(k)(X^\lambda(t_{k-1}))_{\lambda\in\N^d_n\setminus\{0\}}+N(t_k),
\end{equation}
where $A(k) := F(t_k) F(t_{k-1})^{-1}$ and $a(k) := F(t_k) \int_{t_{k-1}}^{t_k} F(s)^{-1}a^c(s) \dd s$. Moreover, if the matrices $(A^c(t))_{t \in \R_+}$ commute, this simplifies to $\smash{a(k) = \int_{t_{k-1}}^{t_k} \exp\big(\int_s^{t_k} A^c(r) \dd r\big) a^c(s) \dd s}$ and $A(k) = \exp\big(\int_{t_{k-1}}^{t_k} A^c(s) \dd s\big)$. In particular, if $X$ is time-homogeneous, $A(k) = \e^{A^c(t_k - t_{k-1})}$ and $\smash{a(k) = \big(\int_0^{t_k - t_{k-1}} \e^{A^c t} \dd t\big) a^c = (A(k) - \mathrm{I})(A^c)^{-1} a^c}$, \mbox{where the latter holds if $A^c$ is regular.}
\item Let $B^c_n(t)$ denote the matrix of coefficients for $X$ from Definition \ref{d:pp} and $F(t)$ the unique solution to the ODE $\mathrm{d}F(t) = B^c_n(t)F(t) \dd t$ with $F(0) = \mathrm{I}$. Then $F(t)$ is regular and
\begin{equation}\label{e:altpssm2}
\E\bigl(X(t_k)^\lambda\big|\F_{t_{k-1}}\bigr) = \smash{\sum_{\mu\in\N^d_{|\lambda|}}} b_{\lambda,\mu}(k)X(t_{k-1})^\mu
\end{equation}

for any $\lambda\in\N^d_n$, where $B_n(k) = (b_{\lambda, \mu}(k))_{\lambda, \mu \in \N^d_n}$ is defined by $B_n(k) := F(t_k)F(t_{k-1})^{-1}$. Moreover, if the matrices $(B^c_n(t))_{t \in \R_+}$ commute, this simplifies to $\smash{B_n(k) = \exp\big(\int_{t_{k-1}}^{t_k} B^c_n(s) \dd s\big)}$. In particular, $B_n(k) = \e^{B_n^c(t_k - t_{k-1})}$ if $X$ is a time-homogeneous polynomial process of order $n$.
\end{enumerate}
\end{proposition}

\begin{proof}
For ease of notation, let $n=1$ so that $(X^\lambda(t))_{\lambda\in\N^d_n\setminus\{0\}} = X(t)$. The local boundedness of the coefficients and the local integrability of $\E(\lVert X(t)\rVert)$ allow us to apply Fubini's theorem in order to obtain that $\mu(t):= \E(X(t) |\F_s)$ is the unique solution to the linear ODE $\mu'(t) = a^c(t) + A^c(t) \mu(t)$ for $t \geq s$ with initial condition $\mu(s) = X(s)$, where uniqueness follows by standard uniqueness results for linear ODEs with locally integrable coefficients. Likewise we have that $\nu(t) := \E((1, X) |\F_s)$ solves $\nu'(t) = B^c_n(t) \nu(t)$ for $t \geq s$ and $\nu(s) = (1, X(s))$. The expressions for $a(k)$, $A(k)$ and $B(k)$ now follow from standard linear systems theory by solving these ODEs, which in particular entails that $F(t)$ is regular, $t \in \R_+$, see \textcite[pp. 253\textendash 255]{kallianpur.13}. The simplification $a(k) = (A(k)-I)(A^c)^{-1} a^c$ for regular $A^c$ follows because $t \sum_{j=1}^\infty \frac{1}{j!} (A^ct)^{j-1} = \int_0^t \e^{A^cs} \dd s$. To derive the expression for $N(t_k)$ in part 1, note that \eqref{eq:longdynamics} implies that $Y(t) := X(t) - \mu(t)$ satisfies the stochastic differential equation $\dd Y(t) = A^c(t) Y(t) \dd t + \dd M(t)$. Since $\mathrm{d}F(t)^{-1} = - F(t)^{-1} A^c(t) \dd t$, we have
\[\mathrm{d}(F^{-1}(t)Y(t)) = F(t)^{-1} \dd Y(t) + \dd F(t)^{-1} Y(t-) = F(t)^{-1} A^c(t) (\Delta Y)(t) \dd t + F(t)^{-1} \dd M(t). \]
By \textcite[Proposition I.1.32]{Jacod2003}, $\{t: \Delta Y(t) \neq 0\}$ is countable for any $\omega \in \Omega$. Hence the first term on the right-hand side above vanishes. The claim about $N$ follows because
\begin{equation}\label{eq: X_closedform}
X(t) = F(t)F(s)^{-1} X(s) + F(t) \int_s^t F(r)^{-1} a^c(r) \dd r + F(t) \int_s^t F(r)^{-1} \dd M(r).
\end{equation}
\end{proof}

\section{Gaussian equivalents of polynomial processes}\label{Xsec5-3}\label{sec:Gauss}

This section is devoted to the construction of Gaussian processes whose mean and covariance function are matched with a given polynomial process. They form a subclass of polynomial models that is well-known to allow for closed-form solutions in filtering problems, linear-quadratic control, and parameter estimation. We start with the discrete-time case.

\subsection{Equivalent Gaussian state space models}\label{Xsec6-3.1}\label{su: GSSM}

\begin{definition}\label{d:gou1}
We call $X$ a \emph{linear Gaussian state space model} if $X(0)$ is Gaussian and
 \begin{equation}\label{e:gou1}
X(t)=a(t)+A(t)X(t-1)+B(t)W(t)
\end{equation}
for $t \in \N^*$, where $W(t)\in\R^d$ is an $\R^d$-valued standard Gaussian random variable that is independent of $\F_{t-1}$ and $a(t)\in\R^d, A(t)\in\R^{d\times d}, B(t)\in\R^{d\times d}$ are deterministic for any $t\in\N^*$.
\end{definition}
Evidently, linear Gaussian state space models are special cases of polynomial state space models in the sense of Definition \ref{d:pssm}, where the noise sequence $N(t)$ is Gaussian. As noted earlier, the key idea of this paper is to replace polynomial models with Gaussian ones sharing the first two moments. To this end, we derive recursions for the second moments.

\begin{lemma}\label{l:recurrence_for_polynomial_process}
Let $X$ be a polynomial state space model of order 1 as in (\ref{e:definition_time_discrete_poly_ssm}) with finite second moments. Then $P(t) := \E(X(t) X(t)^\top)$ is the unique solution to the Lyapunov recurrence relation
\begin{equation}\label{e:covnt}
P(t) =  a(t)a(t)^\top   +   a(t) \mu(t-1)^\top A(t)^\top +    A(t) \mu(t-1) a(t)^\top + A(t)P(t-1)A(t)^\top   +  \Cov(N(t))
\end{equation}
with $P(0) = \E(X(0)X(0)^\top)$, where $\mu(t) = \E(X(t))$. For $P(s, t) := \E(X(t)X(s)^\top)$, $s \leq t$, we have
\[P(s, t) = \Big(\prod_{r=s+1}^{t} A(r)\Big) P(s) + \sum_{r=s+1}^{t} \Big(\prod_{u=r+1}^{t} A(u)\Big)a(r)\mu(s)^\top.\]
\end{lemma}

\begin{proof}
For the first claim, observe that simple algebraic manipulations yield the recurrence
\begin{align}\label{e:algebraic}
X(t)X(t)^\top =~~&a(t)a(t)^\top +  a(t) X(t-1)^\top A(t)^\top + A(t)X(t-1) a(t)^\top  \nonumber\\
& +  A(t)X(t-1) X(t-1)^\top A(t)^\top + N(t) N(t)^\top + a(t) N(t)^\top \nonumber\\
& + N(t) a(t)^\top + N(t) X(t-1)^\top A(t)^\top + A(t)X(t-1) N(t)^\top.
\end{align}
Since $\E(N(t)|\F_{t-1}) = 0$ holds for any $t$, \eqref{e:covnt} follows immediately by taking the expectation on both sides of (\ref{e:algebraic}). Uniqueness follows quickly by considering the difference between two solutions $P_1$ and $P_2$, getting $P_1(t) - P_2(t)=0$ for any $t\in\N$. The second assertion follows from $\E\bigl(X(t)X(s)^\top\bigr)= \E\bigl( \E(X(t)|\F_{s})X(s)^\top\bigr)$ together with the moment formula \eqref{e:moments}.
\end{proof}

Lemma \ref{l:recurrence_for_polynomial_process} allows us to construct a Gaussian counterpart to a given polynomial model.
\begin{proposition}\label{p:representation_discrete_time}\label{p:representation1}
Let $X$ be a polynomial state space model of order 1 as in (\ref{e:definition_time_discrete_poly_ssm})
with finite second moments.
Moreover, consider a linear Gaussian state space model $Y$ as in Definition \ref{d:gou1} with $\E(Y(0))=\E(X(0))$, $\Cov(Y(0))=\Cov(X(0))$, and
coefficients $a(t)$, $A(t)$, $B(t)$ satisfying
\[B(t)B(t)^\top=\Cov(N(t))=:C(t).\]
Then the first two moments of $X$ and $Y$ match, i.e. $\E(X(t))=\E(Y(t))$ and
\[\E\bigl(X(s)X(t)^\top\bigr)=\E\bigl(Y(s)Y(t)^\top\bigr), \quad s \leq t \in \N.\]
\end{proposition}

\begin{proof}
Due to the fact that $\Cov(N(t))$ is nonnegative definite for any $t\in\N$, we can find matrices $B(t)$ satisfying $B(t)B(t)^\top = \Cov(N(t))$. Since $\Cov(B(t)W(t)) = B(t)B(t)^\top = \Cov(N(t))$, the processes $X$ and $Y$ give rise to the same recurrence relations in \eqref{e:moments} and Lemma \ref{l:recurrence_for_polynomial_process}, which characterise the first and second moments.
\end{proof}
We call $Y$ the \emph{Gaussian equivalent}\footnote{We assume $(\Omega, \F, (\F_t)_{t \in \N})$ to be rich enough to carry the Gaussian sequence $(W(t))_{t \in \N^*}$.} of the polynomial state space model $X$. Obviously, it behaves like $X$ in all applications that depend only on the first two moments of the process.

\begin{remark}    If one starts from an arbitrary adapted sequence $X$ with finite second moments, one may wonder whether the existence of a Gaussian equivalent $Y$ sharing the first two moments with $X$ as in Proposition \ref{p:representation_discrete_time} already implies that $X$ is a polynomial state space model in the sense of Definition \ref{d:poly1}. This is not the case: A reconsideration of the proof of Lemma \ref{l:recurrence_for_polynomial_process} shows that if $Y$ is a linear Gaussian state space model as in Proposition \ref{p:representation_discrete_time}, then the first two moments of $X$ and $Y$ match if and only if $X$ is of the form \eqref{e:definition_time_discrete_poly_ssm} with $\E(N(t) X(s)) = 0$ for all $s < t \in \N$ but not necessarily $\E(N(t) \mid \F_{t-1}) = 0$. Therefore the approach of this paper applies also to some non-polynomial processes sharing their moments with a linear Gaussian state space model. However, polynomial state space models are particularly natural candidates for finding such a Gaussian equivalent because it is guaranteed to exist and it can be computed easily in practice by using the polynomial properties.
\label{Xenun5-5}
\end{remark}

\begin{remark}\label{r:covnt}
How do we obtain the covariance matrix $C(t):= \Cov(N(t))$ that is needed to set up the Gaussian equivalent in the previous proposition? It can be determined recursively if $X$ is actually polynomial of order 2. Indeed, note that both the first moments $\mu(t)=\E(X(t))$ as well as the second moments $P(t):= \E(X(t)X(t)^\top)$ can be computed with Proposition \ref{p:moments} or Remark \ref{rem: time-hom} for $n=2$.
Eq. \eqref{e:covnt} now yields the desired covariance matrix
\begin{align*}
\Cov(N(t)) &= P(t) - a(t)a(t)^\top - a(t) \mu(t-1)^\top A(t)^\top - A(t) \mu(t-1) a(t)^\top - A(t)P(t-1)A(t)^\top.
\end{align*}
\end{remark}

\subsection{Equivalent Gaussian Ornstein\textendash Uhlenbeck processes}\label{Xsec7-3.2}

We turn to the continuous-time counterparts of the notions of Section \ref{su: GSSM}. Similarly as before, we fix a filtered probability space $(\Omega, \F, (\F_t)_{t \in \R_+}, \PP)$ with a right-continuous filtration $(\F_t)_{t \in \R_+}$ and an $\mathbb{R}^d$-valued semimartingale $X = (X(t))_{t \in \mathbb{R}_+}$ on this space.

\begin{definition}[Gaussian Ornstein\textendash Uhlenbeck Process]\label{d:gou-cont}
    Let $a^c(t) \in \mathbb{R}^d$, $A^c(t) \in \mathbb{R}^{d \times d}$, and $B(t) \in \R^{d \times d}$ be measurable, deterministic functions such that $\int_0^t \lVert a^c(s)\rVert \dd s < \infty$, $\int_0^t \lVert A^c(s)\rVert \dd s < \infty$ and $\int_0^t \lVert B(s)B(s)^\top \rVert \dd s < \infty$ for any $t \in \R_+$. If $X$ is the (necessarily continuous and unique) solution of the stochastic differential equation
    \begin{equation}
        \mathrm{d}X(t) = \left(a^c(t) + A^c(t)X(t)\right)\dd t + B(t)\dd W(t), \label{eq: GOU}
    \end{equation}
    with $X(0)$ being a Gaussian random variable and $W$ being some $\mathbb{R}^d$-valued standard Wiener process independent of $X(0)$, then $X$ is called a \emph{Gaussian Ornstein\textendash Uhlenbeck process}.
\label{Xenun41-5}
\end{definition}

The following result is the continuous-time version of Lemma \ref{l:recurrence_for_polynomial_process}.

\begin{lemma}\label{lem: p-rec}
    Let $X$ be a polynomial process of order 2 and $M$ be as in \eqref{eq:pp1}. Then there exists a predictable $\R^{d \times d}$-valued process $N^c$ such that $t \mapsto \E(\lVert N^c(t)\rVert)$ is locally bounded and $\mathrm{d}\langle M, M \rangle(t) = N^c(t) \dd t$. Moreover, $P(t) := \mathbb{E}(X(t)X(t)^\top)$ is the unique solution of the Sylvester matrix ODE
        \begin{equation}\label{eq: ODE_P}
            \frac{\mathrm{d}}{\mathrm{d}t}P(t) = \mu(t)a^c(t)^\top + a^c(t)\mu(t)^\top + P(t)A^c(t)^\top + A^c(t)P(t) + \mathbb{E}(N^c(t)), \quad t \in \R_+
        \end{equation}
        with $P(0) = \E(X(0)X(0)^\top)$ and $\mu(t) = \E(X(t))$. Moreover, with initial condition $P(s, s) = P(s)$ for $s \in \R_+$, the function $P(s, t) := \mathbb{E}(X(t)X(s)^\top)$, $t \geq s$, is the unique solution of the linear ODE
        \begin{equation}\label{eq: ODE_Ps}
            \frac{\mathrm{d}P(s, t)}{\dd t}  = \Big(a^c(t) \mu(s)^\top + A^c(t)P(s, t)\Big)
            , \quad t \geq s.
        \end{equation}
\end{lemma}

\begin{proof}
Since $X$ is polynomial of order 2, it admits differential semimartingale characteristics $(b, c, K)$ by the arguments in the proof of \textcite[Theorem 5.2]{Agoitia2020} and it follows from \textcite[Theorem 4.26]{eberlein.kallsen.19} that $N^c(t) = c(t) + \int xx^\top K(t, \mathrm{d}x)$. Moreover, Theorem 5.2 in \textcite{Agoitia2020} implies that $N^c(t)$ is a quadratic polynomial in $X(t)$ whose coefficients are locally bounded in time because $B^c(t)$ is locally bounded. Since $t \mapsto \E(\lVert X(t)\rVert^2)$ is locally bounded by Remark \ref{r:pp}, it follows that $\E(\lVert N^c(t)\rVert)$ is locally bounded as well. Let $i, j \in \{1, \dots, d\}$ and denote by $A^c_{j, \cdot}$ the $j$-th row of $A^c$. Integration by parts yields
        \begin{align*}
        \mathrm{d}(X_iX_j(t)) ={} & X_i(t-)\left(a^c_j(t) + A^c_{j, \cdot}(t)X(t)\right) \dd t + X_j(t-)\left(a^c_i(t) + A^c_{i, \cdot}(t)X(t)\right) \dd t \\
        &{}+ X_i(t-)\dd M_j(t) + X_j(t-) \dd M_i(t) + \mathrm{d}[M_i,M_j](t).
        \end{align*}
        We have $\mathbb{E}\big((\Delta X_i(t))^2\big) = \mathbb{E}(\Delta \langle X_i, X_i \rangle(t)) = \mathbb{E}(\Delta \langle M_i, M_i \rangle(t)) = 0$ and therefore $\Delta X_i(t) = 0$ almost surely for any fixed $t \in \mathbb{R}_+$. By virtue of Fubini's theorem we obtain that
        \begin{align*}
        \mathbb{E}(X_i(t)X_j(t)) ={} & \mathbb{E}(X_i(0)X_j(0)) + \int_0^t \mathbb{E}(N^c(s)_{ij}) \dd s + \int_0^t \left(\mu_i(s)a^c_j(s) + \mu_j(s)a^c_i(s)\right) \dd s \\
        &+{} \int_0^t \left(A^c_{j, \cdot}(s)\mathbb{E}(X_i(s)X(s)) + A^c_{i, \cdot}(s)\mathbb{E}(X_j(s)X(s))\right) \dd s.
        \end{align*}
        Putting together all matrix entries shows that the function $P$ indeed solves the given matrix-valued ODE. Uniqueness is warranted by the global version of the Picard\textendash Lindel\"of theorem. For the proof of the final assertion, we only need to observe that
        \begin{align*}
        \mathbb{E}(X(t)X(s)^\top) = \mathbb{E}(X(s)X(s)^\top) + \int_s^t \left(a^c(u)\E(X(s))^\top + A^c(u)\mathbb{E}(X(u)X(s)^\top)\right) \dd u
        \end{align*}
         for $t \geq s$.
\end{proof}

Finally, we are able to use Lemma \ref{lem: p-rec} to construct the Gaussian Ornstein\textendash Uhlenbeck process sharing the first two moments of our polynomial process under consideration. As before, we assume $(\Omega, \F, (\F_t)_{t \in {\R_+}})$ to be rich enough to carry an $\R^d$-valued standard Wiener process $W$. The following result is the continuous-time counterpart of Proposition \ref{p:representation_discrete_time}.

\begin{proposition}\label{prop: gauss_equiv}
    Let $X$ denote a polynomial process of order 2 as in \eqref{eq:pp1} and let $N^c(t)$ be as in Lemma \ref{lem: p-rec}. Moreover, consider a Gaussian Ornstein\textendash Uhlenbeck process $Y$ as in Definition \ref{d:gou-cont} with $\mathbb{E}(Y(0)) = \mathbb{E}(X(0))$, $\mathrm{Cov}(Y(0)) = \mathrm{Cov}(X(0))$, and coefficients $a^c(t)$, $A^c(t)$, $B(t)$ satisfying
    \[
    B(t)B(t)^\top = C(t) := \mathbb{E}(N^c(t)).
    \]
    Then the first two moments of $X$ and $Y$ coincide, i.e. $\mathbb{E}(X(t)) = \mathbb{E}(Y(t))$ and
    \[
    \mathbb{E}(X(s)X(t)^\top) = \mathbb{E}(Y(s)Y(t)^\top), \quad s \leq t \in \R_+.
    \]
\end{proposition}

\begin{proof}
    By Proposition \ref{l:restriction}, $\E(X(t)) = \E(Y(t))$ both equal $F(t)\int_0^t F(r)^{-1} a^c(r) \dd r + F(t)\E(X(0))$, where $F$ solves the ODE $\mathrm{d}F(t)=A^c(t) F(t) \dd t$ with $F(0) = I$.
    Due to the fact that $C(t)$ is positive semi-definite, we can then find matrices $B(t)$ satisfying $B(t)B(t)^\top = C(t)$ for any $t \in \mathbb{R}_+$. Since $\int_0^t \lVert C(s)\rVert \dd s < \infty$ by Lemma \ref{lem: p-rec}, the same arguments as in the derivation of \eqref{eq: X_closedform} yield that there exists a unique solution $Y$ to \eqref{eq: GOU} with given $Y(0)$. The predictable covariation of the process $\mathrm{d}\widetilde{M}(t) := B(t)\dd W(t)$ satisfies
    \[
    \mathrm{d}\langle \widetilde{M}, \widetilde{M}^\top \rangle (t) = B(t)B(t)^\top \dd t = C(t) \dd t.
    \]
    Consequently, the processes $X$ and $Y$ give rise to the same ODEs \eqref{eq: ODE_P} and \eqref{eq: ODE_Ps} in Lemma \ref{lem: p-rec}, which characterise the second order moments of these processes.
\end{proof}

Parallel to Section \ref{su: GSSM}, we call $Y$ the \emph{Gaussian equivalent} of the polynomial process $X$.

\begin{remark}    The term $C(t) = \E(N^c(t))$ needed for the Gaussian equivalent in Proposition \ref{prop: gauss_equiv} can be obtained easily: Combining the proof of Lemma \ref{lem: p-rec} and Theorem 5.2 in \textcite{Agoitia2020} and equation (6.52) in \textcite{eberlein.kallsen.19}, we obtain
    \begin{equation}\label{eq: ECt}
    C(t)_{i, j} = \sum_{\mu \in \N^d_2} a_{\bm{1}_i + \bm{1}_j, \: \mu}(t) \E(X(t)^\mu),
    \end{equation}
    where $a_{\lambda, \mu}(t) = \sum_{\nu \leq \lambda, \nu \leq \mu} (-1)^{|\nu|} \binom{\lambda}{\nu} b^c_{\lambda - \nu, \mu - \nu}(t)$ with $b^c_{\lambda, \mu}(t)$ being the coefficients for $X$ as in Definition \ref{d:pp}. The terms $\E(X(t)^\mu)$ in \eqref{eq: ECt} \mbox{can be computed using the formulas \eqref{eq:poly_moments_continuous} or \eqref{e:altpssm2}.}
\label{Xenun7-7}
\end{remark}

\section{Filtering, smoothing, and prediction}\label{Xsec8-4}\label{s:filter}

This section is devoted to optimal linear filtering, prediction and smoothing of partially observed polynomial processes. We let either $I := \N$ or $I := \R_+$ and fix a probability space $(\Omega,\F,(\F_t)_{t\in I}, \PP)$ as well as an $\R^d$-valued adapted process $X=(X(t))_{t\in I}$. If $I = \R_+$, we assume $(\F_t)_{t\in I}$ to be right-continuous. Suppose that the components $X_{m+1}(t),\dots,X_d(t)$ are observable whereas $X_1(t),\dots,X_m(t)$ are not. We let the subscript $\oo$ stand for the observable part of a vector $x \in \R^d$ and let $H := (\delta_{m+i, j})_{i=1, \dots, d-m; \: j=1, \dots, d}$, i.e.\ $x_\oo := Hx = (x_{m+1}, \dots, x_d)$. For $\Sigma\in \R^{d\times d}$ we set $\Sigma_{:, \oo} := \Sigma H^\top =\Sigma_{1:d, \: m+1:d}$, $\Sigma_{\oo, :} := H \Sigma =\Sigma_{m+1:d, \: 1:d}$ as well as $\Sigma_{\oo} := H \Sigma H^\top = \Sigma_{m+1:d, \: m+1:d}.$ The subscript u standing for the unobservable part of a vector is treated in the same manner.

We suppose that $\E(\lVert X(t) \rVert^2) < \infty$ for $t \in I$ and consider the following general filtering problem for fixed $t\in I$. The goal is to minimise the mean squared error  $\E(\lVert X(t)-Y\rVert^2)$ over all random variables $Y$ that are measurable with respect to the observable information
\begin{equation}\label{e:info}
\mathscr G_t:= \sigma\big(\big\lbrace X_\oo(s):s\in I, \: s \leq t\big\rbrace\big).
\end{equation}
We call the minimiser of the above minimisation problem the optimal filter for $X$. Regardless of any specific model it is  given by the conditional mean $\smash{\widehat X(t,t) = \E(X(t) |\mathscr G_t).}$

\subsection{Discrete-time linear filtering problems}\label{Xsec9-4.1}\label{su: Kal}

Let $I = \N$. For Gaussian state space models, the \emph{optimal filter} can be computed recursively:
\begin{proposition}[K\'alm\'an filter]\label{p:Ka-lman}
Suppose that $X$ is a linear Gaussian state space model as in Definition \ref{d:gou1} and set $C(t):= B(t)B(t)^\top$. Let $\widehat X(0,-1):= \E(X(0))$, $\widehat\Sigma(0,-1):= \mathrm{Cov}(X(0))$ and
\begin{align*}
\widehat X(t,t)&:= \widehat X(t,t-1)+\widehat\Sigma_{:, \oo}(t,t-1)\widehat\Sigma_{\oo}(t,t-1)^{+}\bigl(X_{\oo}(t)-\widehat X_{\oo}(t,t-1)\bigr),\\
\widehat X(t + 1, t)&:= a(t + 1) + A(t + 1)\widehat X(t, t),  \\
\widehat\Sigma(t,t)&:= \widehat\Sigma(t,t-1)-\widehat\Sigma_{:, \oo}(t,t-1)\widehat\Sigma_{\oo}(t,t-1)^{+}\widehat\Sigma_{\oo, :}(t,t-1),  \\
\widehat\Sigma(t + 1, t)&:= A(t + 1)\widehat\Sigma(t,t)A(t + 1)^\top+C(t + 1)
\end{align*}
for $t \in \N$, where\,$+$\,stands for the Moore\textendash Penrose pseudoinverse of a matrix. Then we have
\begin{align*}
\widehat X(t,t)=\binom{\E(X_{\uu}(t)|\mathscr G_t)}{ X_{\oo}(t)},\qquad
\widehat\Sigma(t,t)=
\begin{pmatrix}
\widehat\Sigma_\uu(t,t) & 0 \\ 0 & 0\end{pmatrix},
\end{align*}
as well as $\E\bigl(\bigl(X(t)-\widehat X(t,t)\bigr)\smash{\bigl(X(t)-\widehat X(t,t)\bigr)}^\top\bigr)=\widehat\Sigma(t,t)$ for $t \in \N$.
\label{Xenun25-6}
\end{proposition}
\begin{proof}
This follows from \textcite[Corollary 6.1]{shumway.stoffer.17} if we choose $x_t$ and $y_t$ in \cite{shumway.stoffer.17} as $X(t)$ and $X_\oo(t)$, respectively. The variables $\Phi_t$, $\Upsilon_t$, $u_t$, $A_t$, $\Gamma_t$, $Q_t$ and $R_t$ in \textcite{shumway.stoffer.17} correspond to $A(t)$, $a(t)$, 1, $H$, 0, $C(t)$ and 0 in our setup. \textcite[Theorem 5.2.1]{Anderson1979} implies that the inverses above can be replaced by pseudoinverses.
\end{proof}

Let us briefly consider the problems of smoothing and prediction which extend the above filtering problem slightly. For prediction the goal is to determine the \emph{optimal predictor}
\[\widehat X(t,s):= \E(X(t)\mid\mathscr G_s)\quad \mbox{for}\quad s< t,\]
i.e.\ $X(t)$ is to be estimated based on the observation only up to time $s< t$.
By contrast, the smoothing or interpolation problem concerns the computation of the \emph{optimal smoother}
\[\widehat X(t,s):= \E(X(t)\mid\mathscr G_s)\quad\mbox{for} \quad s> t,\]
i.e.\ the partially observable signal $X(t)$ is to be estimated in hindsight rather than online as in the filtering case. In this case, one obtains the following counterparts of Proposition \ref{p:Ka-lman}.

\begin{proposition}[Prediction]\label{p:prediction}
 Suppose that $X$ is a linear Gaussian state space model as in Definition \ref{d:gou1} and let $s<t \in \N$. In addition to the objects from Proposition \ref{p:Ka-lman},
 define for $t > s$
\begin{align*}
\widehat X(t, s)&:= a(t)+A(t)\widehat X(t-1,s),  \\
\widehat\Sigma(t, s)&:= A(t)\widehat\Sigma(t-1,s)A(t)^\top+C(t).
\end{align*}
Then we have $\widehat X(t,s)=\E(X(t)\mid\mathscr G_s)$  and $\E\bigl(\big(X(t)-\widehat X(t,s)\big)\big(X(t)-\widehat X(t,s)\big)^\top\bigr)=\widehat\Sigma(t,s).$
\label{Xenun26-7}
\end{proposition}

\begin{proof}
This follows once more from \textcite[Property 6.1,  Corollary 6.1]{shumway.stoffer.17}.
\end{proof}

\begin{proposition}[Smoothing]\label{p:smoothing}
 Suppose that $X$ is a linear Gaussian state space model as in Definition \ref{d:gou1} and let $s>t \in \N$. In addition to the objects from Proposition \ref{p:Ka-lman}, define for $t < s$
 \begin{align*}
G(t)&:= \widehat\Sigma(t,t)A(t + 1)^\top\widehat\Sigma(t+1,t)^{+},\\
\widehat X(t,s)&:= \widehat X(t,t)+G(t)\bigl(\widehat X(t+1,s)-\widehat X(t+1,t)\bigr),\\
\widehat\Sigma(t,s)&:= \widehat\Sigma(t,t)+G(t)\bigl(\widehat\Sigma(t+1,s)-\widehat\Sigma(t+1,t)\bigr)G(t)^\top.
\end{align*}
Then we have $\E\bigl(\bigl(X(t)-\widehat X(t,s)\bigr)\bigl(X(t)-\widehat X(t,s)\bigr)^\top\bigr) =\widehat\Sigma(t,s)$ for any $t < s$ as well as
\begin{align*}
\widehat X(t,s)=\binom{\E(X_\uu(t) |\mathscr G_s)}{ X_\oo(t)},\qquad
\widehat\Sigma(t,s)=
\begin{pmatrix}
\widehat\Sigma_\uu(t,s) & 0 \\ 0 & 0\end{pmatrix}.
\end{align*}
\end{proposition}

\begin{proof}
This is the time-inhomogeneous version of \textcite[Property 6.2]{shumway.stoffer.17}.
\end{proof}

As it is well known, the filtering problem hardly ever allows for a simple solution if one leaves the narrow framework of linear Gaussian models. A common way out is to consider only \emph{linear} filters. That is to say, we try to minimise the mean squared error $\E(\lVert X(t)-Y \rVert^2)$
over all random variables $Y$ that are affine-linear in the observations, i.e.\ of the form
\[Y=\alpha+\sum_{s=0}^t\gamma(s) X_\oo(s)\]
for deterministic $\alpha\in\R^d$ and $\gamma(s)\in\R^{d \times (d-m)}$, $s \leq t$. But even within this class, it is generally not obvious how to determine the corresponding \emph{optimal linear filter} efficiently. In our setup of polynomial models, however, Proposition \ref{p:representation1} opens the door to using Proposition \ref{p:Ka-lman} in the non-Gaussian case: Clearly, in linear Gaussian state space models the optimal filter from Proposition \ref{p:Ka-lman} coincides with the optimal linear filter. While the general optimal filtering problem depends on the whole conditional distribution of $X$ given $X_\oo$,\space{} the objective function of the linear filtering problem depends on the process $X$ only through its first and second order moments of the form $\E(X(t))$ or $\E(X(t) X(s)^\top)$. Thus, Proposition \ref{p:representation1} implies that polynomial state space models admit the same optimal linear filters as their Gaussian equivalent:
\begin{corollary}[Filtering]\label{co:filter}
Let $X$ be a polynomial state space model with finite second moments as in Proposition \ref{p:representation1}.
Define $\widehat X(t,t)$ and $\widehat \Sigma(t, t)$ as in Proposition \ref{p:Ka-lman}. Then $\widehat X(t,t)$ is an optimal linear filter for $X(t)$ and $\E\bigl(\bigl(X(t)-\widehat X(t,t)\bigr)\smash{\bigl(X(t)-\widehat X(t,t)\bigr)}^\top\bigr)=\widehat\Sigma(t,t)$ for $t \in \N$.
\label{Xenun31-1}
\end{corollary}

\begin{proof}
This follows from Proposition \ref{p:Ka-lman} together with Proposition \ref{p:representation1}.
\end{proof}

For prediction and smoothing the situation is essentially the same. For fixed $s, t \in \N$, we try to obtain the best approximation of $X(t)$ within the class of random variables of the form
\[Y=\alpha+\sum_{r=0}^s\gamma(r) X_\oo(r)\]
for deterministic $\alpha\in\R^d$ and $\gamma(r)\in\R^{d \times (d-m)}$, $r \leq s$. We call the corresponding minimiser the \emph{optimal linear predictor} if $s<t$ and the \emph{optimal linear smoother} if $s>t$, respectively.

\begin{corollary}[Prediction]
Let $X$ be a polynomial state space model with finite second moments as in Proposition \ref{p:representation1}. Define $\widehat X(t,s)$ and $\widehat \Sigma(t, s)$ as in Proposition \ref{p:prediction}. Then $\widehat X(t,s)$ is an optimal linear predictor for $X(t)$ and $\E\bigl(\bigl(X(t)-\widehat X(t,s)\bigr)\smash{\bigl(X(t)-\widehat X(t,s)\bigr)}^\top\bigr)=\widehat\Sigma(t,s)$ for $s < t$.
\label{Xenun32-2}
\end{corollary}

\begin{proof}
This follows from Proposition \ref{p:prediction} together with Proposition \ref{p:representation1}.
\end{proof}

\begin{corollary}[Smoothing]\label{co:smooth}
Let $X$ be a polynomial state space model with finite second moments as in Proposition \ref{p:representation1}. Define $\widehat X(t,s)$ and $\widehat \Sigma(t, s)$ as in Proposition \ref{p:smoothing}. Then $\widehat X(t,s)$ is an optimal linear smoother for $X(t)$ and $\E\bigl(\bigl(X(t)-\widehat X(t,s)\bigr)\smash{\bigl(X(t)-\widehat X(t,s)\bigr)}^\top\bigr)=\widehat\Sigma(t,s)$ for $s > t$.
\label{Xenun33-3}
\end{corollary}

\begin{proof}
This follows from Proposition \ref{p:smoothing} together with Proposition \ref{p:representation1}.
\end{proof}

\begin{remark}\label{rem: MSE}
    The above results allow for very simple error estimates for the optimal linear filter of a polynomial state space model. In particular, the mean squared error (MSE) for the optimal linear filter is obtained via $\E(\lVert X(t) - \widehat X(t, s)\rVert^2) = \mathrm{tr}\: \widehat \Sigma(t, s)$.
\end{remark}

As an example we consider the discretely sampled Heston model from Example \ref{ex:heston1}.
\begin{example}[Heston model]\label{ex:heston3}
For ease of exposition we confine ourselves to the case $\delta=\frac{1}{2}$ in (\ref{e:heston1}, \ref{e:heston2}). By \textcite[Theorem 6.25]{eberlein.kallsen.19}, $\widebar X:= (v,Y,Y^2)$ is a time-homogeneous polynomial process because its differential characteristics $(b,c,K)$ are given by
\begin{align*}
b(t)&=
\begin{pmatrix}
\kappa m - \kappa v(t) \\
\mu \\
v(t) + 2\mu Y(t)
\end{pmatrix},\qquad
c(t)=
\begin{pmatrix}
\sigma^2v(t) & \rho \sigma v(t) & 2\sigma Y(t)v(t) \\
\rho \sigma v(t) & v(t) & 2Y(t)v(t) \\
2\sigma Y(t)v(t) & 2Y(t)v(t) & 4Y(t)^2v(t)
\end{pmatrix}, \qquad
K(t) = 0,
\end{align*}
and hence affine-linear respectively quadratic in $\widebar X(t)$. Its restriction $(\widebar X(t_k))_{k \in \N}$ to discrete time with $t_k := k \Delta t$ is a time-homogeneous polynomial state-space model by Proposition \ref{l:restriction}, which implies that $\E(\widebar X(t_k)^\lambda |\F_{t_{k-1}})$ is an affine function of $\widebar X(t_{k-1})^\mu$ with $|\mu|\leq|\lambda|$. Let $X(k):= (v(t_k),\Delta Y(t_k),(\Delta Y(t_k))^2)$. Since $(\Delta Y(t_k))^m$ is a linear combination of $Y(t_k)^\ell Y(t_{k-1})^{m-\ell}$ with $\ell\leq m$, we easily see that $\E(X(k)^\lambda |\F_{t_{k-1}})$ is an affine function of $\widebar X(k-1)^\mu$, $|\mu|\leq|\lambda|$. But since the conditional law of $(v(t_k),\Delta Y(t_k))$ given $\F_{t_{k-1}}$ does not depend on $Y(t_{k-1})$, we actually have that $\E(X(k)^\lambda \mid\F_{t_{k-1}})$ is a function of only $v(t_{k-1})^\mu$ with $\mu\leq|\lambda|$. Hence, $X(k)=(v(t_k),\Delta Y(t_k),(\Delta Y(t_k))^2) $ is a time-homogeneous polynomial state space model as well.

Using the moment formula \cite[Theorem 6.26]{eberlein.kallsen.19} we can determine the (constant) coefficients $a$, $A$ and the matrix $C(k)$ from Proposition \ref{p:representation_discrete_time} for the polynomial state space model $X(k)=(v(t_k),\Delta Y(t_k),(\Delta Y(t_k))^2) $. A straightforward but slightly tedious calculation yields
\begin{align*}
a=
m\begin{pmatrix}
1 - \e^{-\kappa \Delta t} \\ 0 \\ \Delta t - \frac{1}{\kappa}(1 - \e^{-\kappa\Delta t})
\end{pmatrix}, \qquad
A=
\begin{pmatrix}
\e^{-\kappa\Delta t} & 0 & 0  \\
0 & 0 & 0  \\
\frac{1}{\kappa}(1 - \e^{-\kappa\Delta t}) & 0 & 0
\end{pmatrix},
\end{align*}
where we set $\mu = 0$ for convenience. Moreover, $C(k)=(C_{i, j}(k))_{i, j \in \{1, 2, 3\}}$ with

\begin{fitalign}
    {C_{11}(k)} &{\;=\; \frac{(1 - \e^{-\kappa\Delta t}) \sigma^2}{\kappa} \Big[\big(1 - \e^{-\kappa t_k}\big) m + \e^{-\kappa t_k} \mu_v - \big(1 - \e^{-\kappa\Delta t}\big) \frac{m}{2}\Big]} \\
    {C_{12}(k)} &{\;=\;  \frac{m}{\kappa} \rho \sigma \Big[1 - \e^{-\kappa\Delta t}\Big] + \rho \sigma \Big[\mu_v - m\Big] (\Delta t)\e^{-\kappa t_k}} \\
    {C_{13}(k)} &{\;=\;  \frac{\sigma^2}{2\kappa^2} \Big[\big(1 + 4 \rho^2 - \e^{-\kappa\Delta t}\big) m - 2 \big(\mu_v - m\big) \e^{-\kappa t_k}\Big]\Big[ 1 - \e^{-\kappa\Delta t}\Big]} \\
    &{\;\;\;+\; \frac{\sigma^2}{\kappa} \,\Big[\big(1 + \kappa \rho^2\Delta t\big) \big(\mu_v - m\big) (\Delta t)\e^{-\kappa t_k} - 2 \rho^2 m (\Delta t)\e^{-\kappa \Delta t} \Big]} \\
    {C_{22}(k)} &{\;=\; m\Delta t + \frac{1}{\kappa}\Big[\mu_v - m\Big] \Big[\e^{-\kappa\Delta t} - 1\Big]\e^{-\kappa t_k}} \\
    {C_{23}(k)} &{\;=\; \frac{3 \rho \sigma}{\kappa^2} \Big[\big(\mu_v - m\big) \e^{-\kappa t_k} - m \e^{-\kappa\Delta t} \Big] \Big[ \e^{\kappa\Delta t} - 1\Big] - \frac{3 \rho \sigma}{\kappa} \Big[\mu_v - m\Big] \e^{-\kappa t_k} + \frac{3 \rho \sigma}{\kappa} m\Delta t} \\
    {C_{33}(k)} &{\;=\; \bigg[\Big( \frac{2}{\kappa^2} \big(\Sigma_v - 2m\mu_v + m^2\big) - \frac{\sigma^2}{\kappa^3} \big(2 \mu_v - m\big)\Big) \e^{-2\kappa(t_{k-1})} - \frac{\sigma^2}{\kappa^3}\big(\mu_v-m\big)\e^{-\kappa (t_{k-1})} - \frac{\sigma^2}{2\kappa^3} m \bigg]\Big[1 - \e^{-\kappa\Delta t}\Big]^2} \\
    &{\;\;\;+\; \frac{2}{\kappa^3} \Big[3 \sigma^2 (1 + 2 \rho^2) + 2 \kappa^2 m \Delta t \Big]\Big[\mu_v - m\Big]\Big[ \e^{\kappa\Delta t} - 1\Big] \e^{-\kappa t_k} - \frac{6 \sigma^2}{\kappa^2}\Big[ 1 + 2\rho^2 + \kappa \rho^2\Delta t\Big] \Big[\mu_v - m\Big] (\Delta t)\e^{-\kappa t_k}} \\
    &{\;\;\;-\; \frac{3\sigma^2}{\kappa^3} m \Big[ 1 + 8 \rho^2\Big] \Big[1 - \e^{-\kappa\Delta t}\Big] + \frac{12 \sigma^2}{\kappa^2}\rho^2 m (\Delta t)\Big[1 + \e^{-\kappa\Delta t}\Big] + \frac{3\sigma^2}{\kappa^2} m \Delta t + 2 m^2 (\Delta t)^2,}
\end{fitalign}

\noindent where $\mu_v := \E(v(0))$ and $\Sigma_v := \mathrm{Cov}(v(0))$. We can now apply Corollary \ref{co:filter} to filter the squared volatility $v(t)$ from
the observed returns $\Delta Y(t_i)$ and squared returns $(\Delta Y(t_i))^2$, $i=1,\dots,m$ in an optimal linear fashion. Note however that if $\delta \neq \frac{1}{2}$ in \eqref{e:heston2}, $X = (v, \Delta Y, (\Delta Y)^2)$ ceases to be a polynomial state space model and one instead has to consider the polynomial state space model $(v, v^2, \Delta Y, (\Delta Y)^2)$. Optimal linear filters based on higher powers of $\Delta Y$ as well as predictors and smoothers can be determined in the same way. A similar filtering approach in the Heston model was pursued in \textcite{Cacace2019}, see in particular Fig. \ref{fig:heston} there. Note that the filter depends on the parameters $\kappa$, $m$, $\sigma$, $\rho$. An approach to parameter estimation in this setup is established in the companion paper \textcite{kallsen.richert.24b}.

The upper panel of Fig. \ref{fig:heston} contains 2000 daily observations of the Heston variance $v(t)$ (i.e. $\Delta t := \frac{1}{250}$ with $t$ measured in years) with parameters $\kappa = 1$, $m = 0.4^2$, $\sigma = 0.3$ and $\rho = -0.5$. Moreover, it displays the corresponding optimal linear filter $\widehat v(t, t)^{(1)}$ based on observations of only $\Delta Y$ as well as the optimal linear filter $\widehat v(t, t)^{(2)}$ based on observations of both $\Delta Y$ and $(\Delta Y)^2$. Additionally, we plot a particle filter based on an Euler discretisation of the Heston SDE, which estimates the optimal non-linear filter $\E(v(t) \mid \mathcal{G}_t)$ using $N = 10\,000$ simulated particles. The initial distribution used for all three filters is the stationary distribution of the Heston variance $v$. In the middle panel of Fig. \ref{fig:heston}, we replot the filters $\widehat v(t, t)^{(1)}$ and $\widehat v(t, t)^{(2)}$ with bands of one (unconditional) root mean squared error (RMSE). This RMSE was obtained from the filter covariance matrices $\widehat \Sigma(t, t)$ (see Remark \ref{rem: MSE}) for the linear filters $\widehat v(t, t)^{(1)}$ and $\widehat v(t, t)^{(2)}$ and is plotted in the lower panel of Fig. \ref{fig:heston} together with the RMSE for the particle filter, which was computed via Monte Carlo simulation. The linear filter $\widehat v(t, t)^{(2)}$ and the particle filter $\widehat v(t, t)^{\mathrm{part.}}$ are in fact almost indistinguishable in the first panel. This explains why their RMSE is almost identical as can be seen from the bottom panel.
\label{Xenun43-2}
\end{example}

\begin{remark}    Even though one could theoretically include higher powers $(\Delta Y)^\ell$ as observations for the filter, it is apparent from the above Example \ref{ex:heston3} and in particular from Fig. \ref{fig:heston} that the optimal linear filter based on only $\Delta Y$ and $(\Delta Y)^2$ already achieves a comparable accuracy compared to the optimal non-linear filter (approximated by the particle filter $\widehat v(t, t)^{\mathrm{part.}}$). In particular, one cannot expect the optimal linear filter to converge to the optimal non-linear filter in any sense by just including higher and higher powers $(\Delta Y)^\ell$ because the optimal non-linear filter depends on the whole \emph{trajectory} of observations in a non-linear fashion.
\label{Xenun9-9}
\end{remark}

In the case of high-frequency observations one is typically facing microstructure noise, i.e. the observed increments do not really correspond to the increments $\Delta Y(t_k)$ of the logarithmic stock price in Example \ref{ex:heston1}. If this noise is not taken into account, the estimate of the squared volatility $v(t)$ will typically be biased. The following example models this situation.

\begin{example}[Heston model with microstructure noise]
Suppose that we actually observe
\begin{equation}\label{e:withnoise}
\widetilde Y(t_k):= Y(t_k)+\varepsilon(k)
\end{equation}
instead of $Y(t_k)$ in Example \ref{ex:heston1}, where $\varepsilon(k)$, $k\in \N$, denotes a sequence of independent $N(0,\tau^2)$-distributed random variables. Unless $\tau$ is very large, it makes sense to construct a filter that relies linearly on the squared increments $(\Delta \widetilde Y(t_k))^2$. In order to apply Corollary \ref{co:filter}, we need a polynomial state space model including  $(\Delta \widetilde Y(t_k))^2$ as one of its components. To this end, recall from Example \ref{ex:heston3} that $\overline{X}(k) = (v(t_k), Y(t_k), Y(t_k)^2)$ is a time-homogeneous polynomial state space model if $t_k = k \Delta t$ and $\delta = \frac{1}{2}$ in \eqref{e:heston2}. We show that the same holds for the process $X(k) = (v(t_k), Y(t_k), Y(t_k)^2,\widetilde Y(t_k), \widetilde Y(t_k)^2,(\Delta \widetilde Y(t_k))^2)$ and we determine the corresponding matrix $B=(b_{\lambda,\mu})_{\lambda,\mu\in\N^6}$ in the sense of Lemma \ref{l:altpssm}.

\pagebreak

\setcounter{footnote}{2}
\begin{figure}[H]
    \centering
    \includegraphics[width=\textwidth]{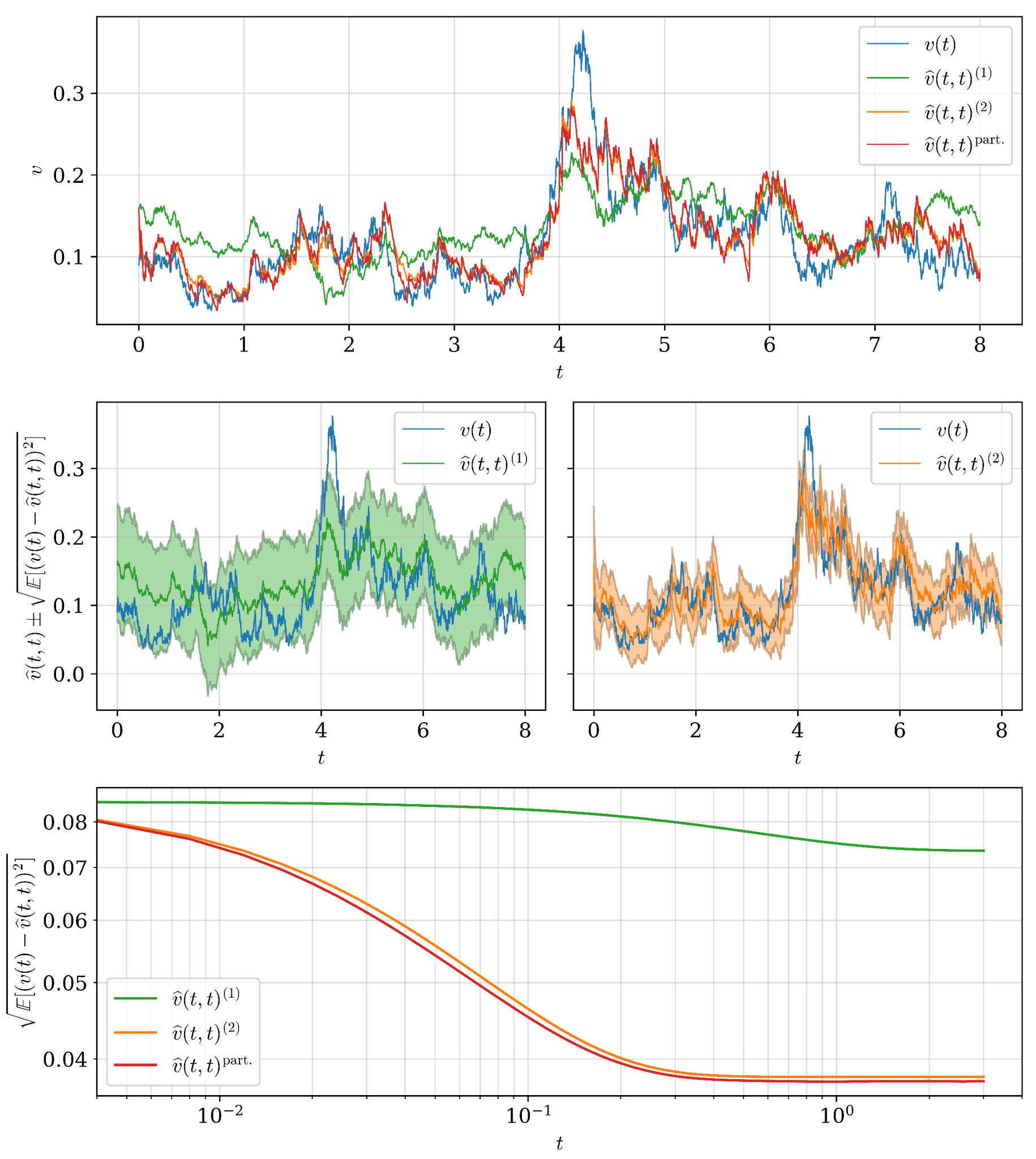}
    \caption{Filters of the Heston variance $v(t)$ and their RMSE. The filter $\widehat v(t, t)^{(1)}$ is the optimal linear filter based on observations of $\Delta Y$  while the filter $\widehat v(t, t)^{(2)}$ is the optimal linear filter based on observations of $\Delta Y$ and $(\Delta Y)^2$. The filter $\widehat v(t, t)^{\mathrm{part.}}$ is a particle filter based on an Euler discretisation of the Heston SDE.\protect\footnotemark}
    \label{fig:heston}
\end{figure}
\footnotetext{Replicating code for this figure can be found in https://github.com/Ambossadore/Quasi-Maximum-Likelihood/blob/master/heston\_filter.py.}

\pagebreak
\noindent \emph{Step 1:}
By \textcite[Theorem 6.25]{eberlein.kallsen.19}, the generating matrix $B^c=(b^c_{\lambda,\mu})_{\lambda,\mu\in\N^3}$ of the time-homogeneous polynomial process $\widebar X(t) = (v(t),Y(t),Y(t)^2)$, $t \in \R_+$, is given by
\[b^c_{\lambda,\mu}=\sum_{\substack{\nu\leq\lambda \\ \nu \leq \mu}}\binom{\lambda}{\nu}a_{\lambda-\nu,\mu-\nu}\]
with
\begin{align*}
a_{\bm{1}_1,0} &= \kappa m, \quad &a_{\bm{1}_1,\bm{1}_1} &= -\kappa, &a_{\bm{1}_2,0} &= \mu, \\
a_{\bm{1}_3,\bm{1}_1} &= 1, \quad &a_{\bm{1}_3,\bm{1}_2} &=2\mu, \\
a_{\bm{1}_{1, 1},\bm{1}_1} &= \sigma^2,\quad
&a_{\bm{1}_{1, 2},\bm{1}_1} &= \rho\sigma,\quad
&a_{\bm{1}_{1, 3}, \bm{1}_{1, 2}} &= \rho\sigma,\\
a_{\bm{1}_{2, 2},\bm{1}_1} &= 1,\quad
&a_{\bm{1}_{2, 3}, \bm{1}_{1, 2}} &= 1,\quad
&a_{\bm{1}_{3, 3},\bm{1}_{1, 3}} &= 1,
\end{align*}
and $a_{\lambda,\mu}=0$ otherwise, where we set $\bm{1}_{i, j} := \bm{1}_i + \bm{1}_j$. Note that $b^c_{\lambda,\mu}=0$ unless $|\lambda| - |\mu| \leq 1$.\\
\\
\emph{Step 2:}
By Proposition \ref{l:restriction}, the matrix $\widebar B=(\widebar b_{\lambda,\mu})_{\lambda,\mu\in\N^3}$ in the sense of Lemma \ref{l:altpssm} of the corresponding time-homogeneous polynomial state space model $(\widebar X(t_k))_{k \in \N}$ is given by $\widebar B= \e^{B^c \Delta t}.$ \\
\\
\emph{Step 3:}
For $\widetilde X(t_k):= (v,Y,Y^2,\widetilde Y, \widetilde Y^2)(t_k)$, $k \in \N$ and $(\lambda,\widetilde\lambda)\in\N^{3+2}=\N^5$
we have that
\begin{align*}
\E(\widetilde X(t_k)^{(\lambda,\widetilde\lambda)}|\F_{t_{k-1}}) =\E\Bigl(v(t_k)^{\lambda_1} Y(t_k)^{\lambda_2 + 2 \lambda_3} \:\E\bigl(\widetilde Y(t_k)^{\widetilde\lambda_1+2\widetilde\lambda_2}\:\big|\:\F_{t_{k-1}}\vee\sigma(v(t_k),Y(t_k))\bigr)\:\Big|\:\F_{t_{k-1}}\Bigr).
\end{align*}
Note that
\begin{align*}
\widetilde Y(t_k)^{\widetilde\lambda_1+2\widetilde\lambda_2} &=\bigl(Y(t_k)+\varepsilon(k)\bigr)^{\widetilde\lambda_1+2\widetilde\lambda_2} =\sum_{\ell=0}^{ \widetilde\lambda_1+2\widetilde\lambda_2}\binom{\widetilde\lambda_1+2\widetilde\lambda_2}{\ell}Y(t_k)^\ell\varepsilon(k)^{\widetilde\lambda_1+2\widetilde\lambda_2-\ell}\\
&=
\begin{cases}
\sum_{\ell=0}^{\widetilde\lambda_1/2+\widetilde\lambda_2}\binom{\widetilde\lambda_1+2\widetilde\lambda_2}{2\ell}(Y(t_k)^2)^\ell \varepsilon(k)^{\widetilde\lambda_1+2\widetilde\lambda_2-2\ell} + \mathrm{odd}
& \text{ if $\widetilde\lambda_1$ is even,}\\
\sum_{\ell=0}^{(\widetilde\lambda_1-1)/2+\widetilde\lambda_2}\binom{\widetilde\lambda_1+2\widetilde\lambda_2}{2\ell + 1}Y(t_k)(Y(t_k)^2)^\ell \varepsilon(k)^{\widetilde\lambda_1-1+2\widetilde\lambda_2-2\ell} + \mathrm{odd}
& \text{ if $\widetilde\lambda_1$ is odd,}
\end{cases}
\end{align*}
where ``odd'' stands for terms involving odd powers of $\varepsilon(k)$.
Using the notation
\begin{align*}
m_p:= \E(\varepsilon(k)^p)=
\begin{cases}
0 & \text{if $p$ is odd}\\
\tau^p\prod_{\ell=1}^{p/2}(2\ell-1)& \text{if $p$ is even,}\\
\end{cases}
\end{align*}
we conclude that $\widetilde X$ is a time-homogeneous polynomial state space model whose corresponding matrix $\widetilde B=(\widetilde b_{\widebar \lambda,\widebar \mu})_{\widebar \lambda, \widebar \mu\in\N^5}$, $\widebar \lambda = (\lambda,\widetilde\lambda)$, $\widebar \mu = (\mu,\widetilde\mu)$, in the sense of Lemma \ref{l:altpssm} is given by
\begin{align*}
\widetilde b_{(\lambda,\widetilde\lambda),(\mu,\widetilde\mu)}
&=\begin{cases}
0 & \text{if }\widetilde\mu\neq0,\\
\sum_{\ell=0}^{\widetilde\lambda_1/2+\widetilde\lambda_2}\binom{\widetilde\lambda_1+2\widetilde\lambda_2}{2\ell}m_{\widetilde\lambda_1+2(\widetilde\lambda_2-\ell)}
\overline b_{\lambda+\ell \bm{1}_3,\mu}
& \text{if $\widetilde\lambda_1$ is even,}\\
\sum_{\ell=0}^{(\widetilde\lambda_1 - 1)/2+\widetilde\lambda_2}\binom{\widetilde\lambda_1+2\widetilde\lambda_2}{2\ell+1}m_{\widetilde\lambda_1-1+2(\widetilde\lambda_2-\ell)}
\overline b_{\lambda+\bm{1}_2+\ell \bm{1}_3,\mu}
& \text{if $\widetilde\lambda_1$ is odd.}
\end{cases}
\end{align*}

\noindent\emph{Step 4:}
Finally, set $X(t_k):= (v,Y,Y^2,\widetilde Y, \widetilde Y^2,(\Delta \widetilde Y)^2)(t_k)$ for $k \in \N$.
Since

\begin{fitalign}
(\Delta \widetilde Y(t_k))^{2m}
=&\sum_{\ell=0}^{2m} (-1)^\ell \binom{2m}{\ell}\widetilde Y(t_k)^\ell \widetilde Y(t_{k-1})^{2m-\ell}\\
=&\sum_{\ell=0}^{m}\binom{2m}{2\ell}(\widetilde Y(t_k)^2)^\ell (\widetilde Y(t_{k-1})^2)^{m-\ell} -\sum_{\ell=0}^{m-1}\binom{2m}{2\ell+1}\widetilde Y(t_k)(\widetilde Y(t_k)^2)^\ell \widetilde Y(t_{k-1})(\widetilde Y(t_{k-1})^2)^{m-1-\ell},
\end{fitalign}

\pagebreak
\noindent a straightforward computation shows that $(X(t_k))_{k \in \N}$ is a time-homogeneous polynomial state space model of any order whose matrix $B=(b_{\lambda,\mu})_{\lambda,\mu\in\N^6}$
\mbox{in the sense of Lemma \ref{l:altpssm} is given by}
\begin{align*}
b_{(\lambda,\overline\lambda),(\mu,\overline\mu)}
&=\begin{cases}
0 & \text{if }\overline\mu\neq0,\\
\sum_{\ell=0}^{\overline\lambda}\binom{2\overline\lambda}{2\ell}\widetilde b_{\lambda+\ell\bm{1}_5,\;\mu+(\overline\lambda-\ell)\bm{1}_5}
- \sum_{\ell=0}^{\overline\lambda-1}\binom{2\overline\lambda}{2\ell+1}\widetilde b_{\lambda+\bm{1}_4+\ell\bm{1}_5, \;\mu+\bm{1}_4+(\overline\lambda-1-\ell)\bm{1}_5}
& \text{otherwise}
\end{cases}
\end{align*}
for $(\lambda,\overline\lambda),(\mu,\overline\mu)\in\N^{5+1}=\N^6$.

\noindent \emph{Step 5:} For the application to filtering it is enough to consider the order $n=2$.
The (constant) coefficients
$a(t)=(a_1,\dots,a_6)$ and
$A(t)=(A_{j\ell})_{j,\ell=1,\dots,6}$ from Proposition \ref{p:representation_discrete_time} for $X$
are given by $a_j= b_{\bm{1}_j,0}$ and $A_{i, j}= b_{\bm{1}_i,\bm{1}_j}$, see Lemma \ref{l:altpssm}.
Moreover, the matrix $C(k)=\Cov(N(t_k))$ can be computed as explained in Remark \ref{r:covnt}. We can then apply Corollary \ref{co:filter} to filter the squared volatility $v(t)$ and, if desired,
also the logarithmic stock price $Y(t)$ from the observations
$\widetilde Y(t_k)$, $\widetilde Y(t_k)^2$, $(\Delta \widetilde Y(t_k))^2$, $k=1,\dots,m$.
\label{Xenun44-3}
\end{example}

\subsection{Continuous-time filtering of polynomial processes}\label{Xsec10-4.2}\label{su: KB}

The ideas of the previous section allow for a natural continuous-time counterpart. Accordingly, let now $I = \R_+$. Analogous to the previous subsection, the filtering problem possesses an explicit solution for Gaussian Ornstein\textendash Uhlenbeck processes:

\begin{proposition}[K\'alm\'an\textendash Bucy filter]\label{prop: KalBuc}
Suppose that $X$ is a Gaussian Ornstein\textendash Uhlenbeck process as in Definition \ref{d:gou-cont} and define $\mu(0) := \E(X(0))$ as well as $\Sigma(0) := \mathrm{Cov}(X(0))$. Set $C(t) := B(t)B(t)^\top$ and assume that $\mathrm{col} \:A^c_{\oo, \uu}(t) \subseteq \mathrm{col} \:C_\oo(t)$ holds for almost every $t \in \R_+$. Assume moreover that $\smash{\int_0^t \lVert (C_\oo(s)^+)^{\frac{1}{2}} a^c_{\oo}(s)\rVert^2 \dd s< \infty}$ and $\smash{\int_0^t \lVert (C_\oo(s)^+)^{\frac{1}{2}} A^c_{\oo, :}(s)\rVert^2 \dd s< \infty}$ for every $t \in \R_+$. Define $\Phi(t) := [A^c(t) - C_{:, \oo}(t) C_\oo(t)^+ A^c_{\oo, :}(t)]$ as well as $R(t) := C(t) - C_{:, \oo}(t) C_\oo(t)^+C_{\oo, :}(t)$. Then the matrix Riccati differential equation
\begin{equation}\label{eq: Kalman_eq_P}
\frac{\mathrm{d}}{\mathrm{d}t}\widehat \Sigma(t) = \Phi(t) \widehat{\Sigma}(t) + \widehat{\Sigma}(t) \Phi(t)^\top - \widehat{\Sigma}(t) A^c_{\oo, :}(t)^\top C_\oo(t)^+ A^c_{\oo, :}(t) \widehat{\Sigma}(t) + R(t)
\end{equation}
with $\widehat{\Sigma}(0) = \Sigma(0) - \Sigma_{:, \oo}(0) \Sigma_\oo(0)^+ \Sigma_{\oo, :}(0)$, where\,$+$\,denotes the Moore\textendash Penrose pseudoinverse, has a unique absolutely continuous, symmetric, positive semidefinite solution $\widehat{\Sigma}: \R_+ \to \R^{d \times d}$. Let $\Psi(t) := \big[\widehat{\Sigma}(t) A^c_{\oo, :}(t)^\top + C_{:, \oo}(t)\big] C_\oo(t)^+$. Set $\widehat X(0) = \mu(0) + \Sigma_{:, \oo}(0) \Sigma_\oo(0)^{+}(X_\oo(0) - \mu_\oo(0))$ and let $\widehat X$ be the unique corresponding $\mathbb{R}^d$-valued solution to the linear stochastic differential equation
\begin{equation}\label{eq: Kalman_eq}
\mathrm{d}\widehat X(t) = \big[a^c(t) - \Psi(t) a^c_\oo(t) + (A^c(t) - \Psi(t) A^c_{\oo, :}(t)) \widehat X(t)\big] \dd t + \Psi(t) \dd X_{\oo}(t).
\end{equation}
Then we have $\E\bigl(\bigl(X(t)-\widehat X(t)\bigr)\smash{\bigl(X(t)-\widehat X(t)\bigr)}^\top\bigr)=\widehat{\Sigma}(t)$ for any $t \in \R_+$ as well as
\begin{align*}
\widehat X(t)=\binom{\E(X_{\uu}(t)|\mathscr G_t)}{ X_{\oo}(t)},\qquad
\widehat{\Sigma}(t)=
\begin{pmatrix}
\widehat{\Sigma}_\uu(t) & 0 \\ 0 & 0\end{pmatrix}.
\end{align*}
\end{proposition}

\begin{proof}
For brevity, we write $L_{\mathrm{loc}}^p$ for the space of measurable functions $h: \R_+ \to \R$ such that $\int_0^t \lVert h(s)\rVert^p \dd s < \infty$ for all $t \in \R_+$. First, the properties of the normal distribution imply the given expressions for the initial values $\widehat X(0) = \E(X(0) |X_\oo(0))$ and $\widehat{\Sigma}(0) = \mathrm{Cov}(X(0)|X_\oo(0))$. Set $D(t):=C_\oo(t)$ and $P(t):=D(t)D(t)^+$, the orthogonal projection onto $\operatorname{col}D(t)$. Since $C(t) \geq 0$, $C_{:,\oo}(t)D(t)^+C_{\oo,:}(t)\leq C(t),$ and consequently $C_{:,\oo}D^{+1/2}\in L^2_{\mathrm{loc}}$. This implies that $\Phi$, $R$ and $(A_{\oo,:}^c)^\top D^+A_{\oo,:}^c$ are elements of $L^1_{\mathrm{loc}}$. Hence the Riccati equation \eqref{eq: Kalman_eq_P} possesses a unique absolutely continuous symmetric positive semidefinite solution \(\widehat\Sigma\). Moreover,
$$ \Psi D^{1/2} = \widehat\Sigma(A_{\oo,:}^c)^\top D^{+1/2} + C_{:,\oo}D^{+1/2} \in L^2_{\mathrm{loc}}, $$
while $\Psi a_\oo^c$ and $\Psi A_{\oo,:}^c$ are elements of $L^1_{\mathrm{loc}}$. Thus, the equation \eqref{eq: Kalman_eq} is well-defined and has a unique corresponding solution.

It remains to verify that this solution is the conditional mean. Since $\operatorname{col}A^c_{\oo,\uu}(t)\subseteq\operatorname{col}C_\oo(t),$ we have $P(t)A^c_{\oo,\uu}(t)=A^c_{\oo,\uu}(t)$ almost everywhere. Moreover, $P(t)B_{\oo,:}(t)=B_{\oo,:}(t)$ because $\operatorname{col}B_{\oo,:}(t)=\operatorname{col}C_\oo(t)$. On an enlargement of the probability space carrying an independent $(d-m)$-dimensional Wiener process $W'$, define the process $\widetilde W$ via
$$ \dd\widetilde W(t) := C_\oo(t)^{+1/2}B_{\oo,:}(t)\dd W(t) +(I-P(t))\dd W'(t). $$
Then $\mathrm{d}\langle\widetilde W,\widetilde W^\top\rangle(t) = \bigl(P(t)+I-P(t)\bigr)\dd t=I\dd t,$ so $\widetilde W$ is a standard Wiener process with respect to $\mathscr{F}_t \vee \sigma(\{U(s): 0 \leq s \leq t\})$, where $U(t) := \int_0^t (I - P(s)) \dd W'(s)$. Furthermore, let $Y(0) := X_\oo(0)$ and $$\dd Y(t):= C_\oo(t)^{+1/2} \bigl[\mathrm{d}X_\oo(t)-(a_\oo^c(t)+A^c_{\oo,\oo}(t)X_\oo(t))\dd t\bigr] +(I-P(t))\dd W'(t),$$ we obtain $\dd Y(t) = C_\oo(t)^{+1/2}A^c_{\oo,\uu}(t)X_\uu(t)\dd t +\dd \widetilde W(t).$ The additional component $\int (I-P)\dd W'$ is independent of $X$, while the remaining component generates the same filtration as $X_\oo$. Hence $\E(X(t) \mid \mathscr{G}_t) = \E(X(t) \mid \sigma(\{Y(s): 0 \leq s \leq t\}))$. The claim then follows by applying the arguments in \textcite[Theorem 10.5.1]{kallianpur.13} or \textcite[Section 4.4]{Davis1977} to the system $(X_\uu(t), Y(t))_{t \in \R_+}$ with $D^+(t)$ in place of $D^{-1}(t)$.
\end{proof}

The optimal predictor and smoother $\widehat{X}(t,s) := \mathbb{E}(X(t) |\mathscr{G}_s)$ can be obtained similarly.

\begin{proposition}[Prediction]\label{prop: KalBuc_pred}
Suppose that $X$ is a Gaussian Ornstein--Uhlenbeck process as in Definition \ref{d:gou-cont} and adopt the notation and assumptions from Proposition \ref{prop: KalBuc}. Fix $s \in \mathbb{R}^+$. Then the equation
\begin{equation}\label{eq: Kalman_pred_P}
\frac{\mathrm{d}}{\mathrm{d}t} \widehat{\Sigma}(t, s) = A^c(t) \widehat{\Sigma}(t, s) + \widehat{\Sigma}(t, s) A^c(t)^\top + C(t)
\end{equation}
for $t \geq s$ with $\widehat{\Sigma}(s, s) = \widehat{\Sigma}(s)$ has a unique absolutely continuous, symmetric, positive semidefinite solution $\widehat{\Sigma}(t, s)$. Together with $\widehat X(s, s) = \widehat X(s)$, let $\widehat X(t, s)$ be the unique $\mathbb{R}^d$-valued solution to
\begin{equation}\label{eq: Kalman_pred}
\mathrm{d}\widehat X(t, s) = \left( a^c(t) + A^c(t) \widehat X(t, s) \right) \dd t
\end{equation}
for $t \geq s$. Then $\widehat X(t, s) = \E(X(t) |\mathscr G_s)$ and $\E\bigl(\big(X(t)-\widehat X(t,s)\big)\big(X(t)-\widehat X(t,s)\big)^\top\bigr)=\widehat{\Sigma}(t,s).$
\label{Xenun29-10}
\end{proposition}

\begin{proof}
    This follows quickly by taking conditional expectations given $\mathscr G_s$ on both sides of the closed-form solution \eqref{eq: X_closedform} with $\mathrm{d}M(t) = B(t) \dd W(t)$ or as in \textcite[p. 142]{Davis1977}.
\end{proof}

\begin{proposition}[Smoothing]\label{prop: KalBuc_smooth}
Suppose that $X$ is a Gaussian Ornstein--Uhlenbeck process as in Definition \ref{d:gou-cont} and adopt the notation and assumptions from Proposition \ref{prop: KalBuc}. Assume moreover that $\widehat \Sigma(t)^+$ is locally bounded and fix $s \in \R_+$. Then the equation
\begin{equation}\label{eq: Kalman_smooth_P}
    \frac{\mathrm{d}}{\mathrm{d}t}\widehat{\Sigma}(t, s) = \big(\widehat \Sigma(t) \Phi(t)^\top + R(t)\big) \widehat{\Sigma}(t)^+ \widehat \Sigma(t, s) + \widehat \Sigma(t, s) \widehat{\Sigma}(t)^+ \big(\Phi(t)\widehat \Sigma(t) + R(t)\big) - R(t)
\end{equation}
for $t \leq s$ with $\widehat \Sigma(s, s) = \widehat \Sigma(s)$ has a unique absolutely continuous, symmetric, positive semidefinite solution $\widehat \Sigma(t, s)$. Moreover, define $\varphi(t) := a^c(t) - C_{:, \oo}(t) C_\oo(t)^{+} a^c_\oo(t)$. Together with the terminal value $\widehat X(s, s) = \widehat X(s)$, let $\widehat X(t, s)$ denote the unique $\R^d$-valued solution to the equation
\begin{equation}\label{eq: Kalman_smooth}
\mathrm{d}\widehat{X}(t, s) = \Big[\varphi(t) + \Phi(t) \widehat X(t, s) + R(t) \widehat{\Sigma}(t)^+ \big(\widehat{X}(t, s) - \widehat{X}(t)\big)\Big] \dd t + C_{:, \oo}(t) C_\oo(t)^+ \dd X_\oo(t)
\end{equation}
for $t \leq s$. Then we have $\E\bigl(\bigl(X(t)-\widehat X(t,s)\bigr)\bigl(X(t)-\widehat X(t,s)\bigr)^\top\bigr) = \widehat{\Sigma}(t,s)$  as well as
\begin{align*}
\widehat X(t,s)=\binom{\E(X_\uu(t) |\mathscr G_s)}{ X_\oo(t)},\qquad
\widehat{\Sigma}(t,s)=
\begin{pmatrix}
\widehat{\Sigma}_\uu(t,s) & 0 \\ 0 & 0\end{pmatrix}.
\end{align*}
\end{proposition}

\begin{proof}
For convenience we may assume that $a^c \equiv 0$ and that $\E(X(0)\mid X_\oo(0))=0$, because the general case can be deduced as a simple corollary. Set $D(t):=C_\oo(t)$ and $P(t):=D(t)D(t)^+$ and let $W'$ be an independent $\R^{d-m}$-valued standard Wiener process on an enlargement of the underlying probability space. Define
\[\dd\widetilde W(t) := D(t)^{+1/2} \big[\dd X_\oo(t)-A^c_{\oo,:}(t)\widehat X(t)\dd t\big] + \big(I-P(t)\big)\dd W'(t).\]
Since $\widehat X_\oo=X_\oo$, the assumption $\mathrm{col}\:A^c_{\oo,\uu}(t)\subseteq\mathrm{col}\:D(t)$ yields $(I-P(t))A^c_{\oo,:}(t)\big(X(t)-\widehat X(t)\big)=0$ for almost every $t$. Moreover, $\mathrm{col}\:B_{\oo,:}(t)=\mathrm{col}\:D(t)$ and hence $P(t)B_{\oo,:}(t)=B_{\oo,:}(t)$. Consequently,
\[\dd\widetilde W(t)=D(t)^{+1/2}A^c_{\oo,:}(t)\widetilde X(t)\dd t+D(t)^{+1/2}B_{\oo,:}(t)\dd W(t)+(I-P(t))\dd W'(t),\]
where $\widetilde X(t):=X(t)-\widehat X(t)$, and $\dd\langle\widetilde W,\widetilde W^\top\rangle(t) = \big(P(t)+I-P(t)\big)\dd t =I\dd t$. Thus $\widetilde W$ is a standard Wiener process with respect to the filtration $\mathscr{G}_t \vee \sigma(\{W'(s): 0 \leq s \leq t\})$. Since $W'$ is independent of $X$, conditioning with respect to this enlarged observation
filtration does not change the conditional expectations of $X$. Therefore, by \textcite[Lemma 6.2.7]{oeksendal.03} (or the multivariate extension thereof), it follows $\widehat X(t, s) = \int_0^s \frac{\mathrm{d}}{\mathrm{d}r} \E\bigl(X(t) \widetilde W(r)^\top\bigr) \dd \widetilde W(r)$. This implies
\begin{equation}\label{e:stern2}
    \widehat X(t, s) = \widehat X(t) +  \int_t^s \frac{\mathrm{d}}{\mathrm{d}r} \E\bigl(X(t) \widetilde W(r)^\top\bigr) \dd \widetilde W(r).
\end{equation}
As in \textcite[Lemma 6.2.6]{oeksendal.03} it follows that 
\begin{fitequation}\label{eq: decomp}
\E\bigl(X(t)\widetilde W(r)^\top\bigr) =\int_0^r \E\bigl(X(t)\widetilde X(u)^\top\bigr) A^c_{\oo,:}(u)^\top D(u)^{+1/2}\dd u + \E\left[X(t)\left(\int_0^rD(u)^{+1/2}B_{\oo,:}(u)\dd W(u)\right)^\top\right].
\end{fitequation}
For $r > t$ the second summand on the right is independent of $r$ because $X(t)$ is independent of $\int_t^r D(u)^{+1/2} B_{\oo, :}(u) \dd W(u)$. Therefore, the derivative of the second summand in \eqref{eq: decomp} with respect to $r$ is 0 for $r > t$. The derivative of the first summand in \eqref{eq: decomp} with respect to $r$ is $\E(X(t) \widetilde X(r)^\top) A^c_{\oo, :}(r)^\top D(r)^{+1/2}$. Let $\widetilde A^c(t) := A^c(t) - \Psi(t) A^c_{\oo, :}(t)$. By \eqref{eq: Kalman_eq}, $\widetilde X(u)$ fulfils the stochastic differential equation
    \[\mathrm{d}\widetilde X(t) = \widetilde A^c(t) \widetilde X(t) \dd t + (B(t) - \Psi(t) B_{\oo, :}(t)) \dd W(t),\]
    whence $\widetilde X(r) = \Gamma(r)^\top \Gamma^{-1}(t)^\top \widetilde X(t) + \Gamma(r)^\top \int_t^r \Gamma^{-1}(u)^\top (B(u) - \Psi(u) B_{\oo, :}(u)) \dd W(u)$, where $\Gamma(t) \in \R^{d \times d}$ solves $\mathrm{d}\Gamma(t) = \Gamma(t)\widetilde A^c(t)^\top \dd t$ with $\Gamma(0) = \mathrm{I}_d$. It follows that $\E(X(t) \widetilde X(r)^\top) = \widehat{\Sigma}(t) \Gamma(t)^{-1} \Gamma(r)$ because $\E(X(t) \widetilde X(t)^\top) = \E(\widetilde X(t) \widetilde X(t)^\top) = \widehat{\Sigma}(t)$. All in all we get
\begin{equation}\label{eq: Oksendal}
\frac{\dd}{\dd r} \E\bigl(X(t)\widetilde W(r)^\top\bigr) = \widehat\Sigma(t)\Gamma(t)^{-1}\Gamma(r) A^c_{\oo,:}(r)^\top D(r)^{+1/2}.
\end{equation}
Inserting Eq. \eqref{eq: Oksendal} into \eqref{e:stern2} yields $\widehat X(t, s) = \widehat X(t) + \widehat{\Sigma}(t) \lambda(t)$, where $\lambda(t)$ solves \[\lambda(t) = \int_t^s \widetilde A^c(r)^\top \lambda(r) \dd r + \int_t^s A^c_{\oo, :}(r)^\top D(r)^{+1/2} \dd \widetilde W(r).\] An application of the product rule and \eqref{eq: Kalman_eq_P} show that
    \begin{align*}
        \widehat \Sigma(t) \lambda(t) &= -\int_t^s \Big( \frac{\mathrm{d}}{\mathrm{d}r}\widehat \Sigma(r) - \widehat \Sigma(r) \widetilde A^c(r)^\top\Big) \lambda(r) \dd r + \int_t^s \widehat \Sigma(r) A^c_{\oo, :}(r)^\top D(r)^{+1/2} \dd \widetilde W(r) \\
        &= -\int_t^s \Big( \Phi(r) \widehat{\Sigma}(r) + R(r)\Big) \lambda(r) \dd r + \int_t^s \widehat \Sigma(r) A^c_{\oo, :}(r)^\top D(r)^{+1/2} \dd \widetilde W(r).
    \end{align*}
    Plugging this and \eqref{eq: Kalman_eq} into $\widehat X(s) - \widehat X(t, s) = \widehat X(s) - \widehat X(t) - \widehat \Sigma(t) \lambda(t)$ yields
    \begin{fitequation*}
        \widehat X(s) - \widehat X(t, s) = \int_t^s \Big[A^c(r) \widehat X(r, s) + R(r) \lambda(r) - C_{:, \oo}(r) D(r)^{+} A^c_{\oo, :}(r) \widehat \Sigma(r) \lambda(r) \Big] \dd r + \int_t^s C_{:, \oo}(r) D(r)^{+1/2} \dd \widetilde W(r).
    \end{fitequation*}
Here, we can replace $\widehat \Sigma(r) \lambda(r)$ by $\widehat X(r, s) - \widehat X(r)$. Moreover, since $\widehat\Sigma^+$ is locally bounded and $\widehat\Sigma$ is continuous, $\operatorname{rank}\widehat\Sigma(t)$ is locally constant. Hence, for almost every $t$ and every $v\in\ker\widehat\Sigma(t)$, we have $v^\top\widehat\Sigma'(t)v=0$. By \eqref{eq: Kalman_eq_P}, this gives $v^\top R(t)v=0.$ Since $R(t)$ is positive semidefinite, $R(t)v=0$ and consequently $R(t) = R(t)\widehat\Sigma(t)^+\widehat\Sigma(t)$ almost everywhere. Therefore $R(t)\lambda(t) = R(t)\widehat\Sigma(t)^+ \big(\widehat X(t,s)-\widehat X(t)\big).$ Together, using $C_{:, \oo}(t) D(t)^{+1/2} \dd \widetilde W(t) = C_{:, \oo}(t) D(t)^+[\mathrm{d}X_\oo(t) - A^c_{\oo, :}(t) \widehat{X}(t)\dd t]$, this yields \eqref{eq: Kalman_smooth}. Now, note that $\widetilde X(t, s) := X(t) - \widehat X(t, s)$ is of the form $\smash{\widetilde X(t, s) = \widetilde X(t) - \int_t^s \widehat{\Sigma}(t) \Gamma(t)^{-1} \Gamma(r) A_{\oo, :}^c(r)^\top D(r)^{+1/2} \dd \widetilde W(r)}$. The identity \eqref{eq: Kalman_smooth_P} follows by applying integration by parts to $\widetilde X(t, s) \widetilde X(t, s)^\top$. In particular, the coefficients of \eqref{eq: Kalman_smooth_P} are locally integrable by the assumptions of Proposition \ref{prop: KalBuc} and the local boundedness of $\widehat\Sigma^+$, so \eqref{eq: Kalman_smooth_P} has a unique solution.
\end{proof}
\begin{remark}    From \eqref{eq: Kalman_smooth} it follows that the optimal smoother $\widehat X_\uu(t, s)$ of a Gaussian Orn\-stein--Uhlenbeck process is absolutely continuous in $t < s$ if and only if $C_{\uu, \oo}(t) = 0$ for almost every $t < s$. Indeed, the rowspace of $C_{\uu, \oo}(t)$ is contained in the rowspace of $C_\oo(t)$ because $C(t)$ is positive semidefinite, and so $C_{\uu, \oo}(t) = C_{\uu,\oo}(t) C_\oo(t)^+ C_\oo(t)$. Now, $\mathrm{tr}[C_{\uu, \oo}(t) C_\oo(t)^+ C_{\oo,\uu}(t)] = \lVert C_{\uu, \oo}(t) (C_\oo(t)^+)^{1/2}\rVert_F^2$, so $C_{\uu, \oo}(t) C_\oo(t)^+ C_{\oo,\uu}(t) = 0$ if and only if $C_{\uu, \oo}(t) = 0$, and thus the local martingale part of $C_{\uu, \oo}(t) C_\oo(t)^+ \dd X_\oo(t)$ vanishes if and only if $C_{\uu, \oo}(t) = 0$.
\label{Xenun10-10}
\end{remark}

In order to tackle filtering problems for polynomial processes beyond the Gaussian case, we focus as in the previous section on minimising the error $\mathbb{E}(\lVert X(t) - Y\rVert^2)$ within the class of all random variables linear in the observations. Here we define ``linear'' in the following generalised sense: Let $X$ denote an $\R^d$-valued process with $\mathbb{E}(\|X(t)\|^2) < \infty$ for any $t \in \mathbb{R}^+$, and fix $d'$. For any $\R^{d'}$-valued random variable $Z$, we let $\|Z\|_{2} := \sqrt{\mathbb{E}(Z^\top Z)}$. We call the $\|\cdot\|_2$-closure of the set
\begin{equation}\label{eq: char_E}
\mathcal{E}(X,t,d') := \left\{ \alpha + \sum_{i=1}^n \gamma_i X(t_i) \;\big|\; n \in \mathbb{N},\: \alpha \in \R^{d'}, \: \gamma_i \in \R^{d' \times d}, 0 \leq t_i \leq t, i = 1, \ldots, n \right\},
\end{equation}
equipped with the norm $\|\cdot\|_2$, the \emph{linear space of $X$ at time $t$}. Accordingly we write
\[
L^2(X, t, d') := \smash{\overline{\mathcal{E}(X, t, d')}}^{\|\cdot\|_2}.
\]
The linear filtering, prediction or smoothing problems therefore naturally entail minimising the mean squared distance $\mathbb{E}(\lVert X(t) - Z\rVert^2)$ over all random variables $Z \in L^2(X_\oo, s, d)$. We introduce the following auxiliary result:

\begin{lemma}\label{lem: funcana} Let $X$ be an $\R^d$-valued polynomial process of order 2 as in Eq. \eqref{eq:pp1}. Furthermore, let $Y$ denote a Gaussian equivalent of $X$ as in Proposition \ref{prop: gauss_equiv} and set $q := d - m$. Then we have:
\begin{enumerate}
    \item Assume that $C_\oo(t)$ is regular for $t \in \R_+$ and $C_\oo(t)^{-1}$ is locally bounded. Then, for $t \in \R_+$,
    \begin{fitequationtwo}\label{eq: char_L2}
    L^2(X_\oo, t, d') = \Big\lbrace \alpha + \beta X_\oo(0) + \int_0^t \gamma(s) \dd X_\oo(s) \: \Big| \: \alpha \in \R^{d'},\: \beta \in \R^{d' \times q}, \: \gamma \in L^2\big([0, t], \R^{d'\times q}\big) \Big\rbrace.
    \end{fitequationtwo}
    \item Fix $t \in \R_+$. Then there exists an isometric isomorphism $I_t : L^2(X, t, d') \to L^2(Y, t, d')$ with $I_t\big( \alpha + \sum_{i=1}^n \gamma_i X(t_i)\big) = \alpha + \sum_{i=1}^n \gamma_i Y(t_i)$ for $0 \leq t_i \leq t$ and $I_t\big(L^2(X_\oo, s, d')\big) = L^2(Y_\oo, s, d')$ for $s \leq t$. If $C_\oo(s)$ is regular for any $s \in \R_+$ and $C_\oo(s)^{-1}$ is locally bounded, then
    \begin{equation}\label{eq: isometry}
        I_t\Big[ \alpha + \beta X_\oo(0) + \int_0^s \gamma(r) \dd X_\oo(r)\Big] = \alpha + \beta Y_\oo(0) + \int_0^s \gamma(r) \dd Y_\oo(r)
    \end{equation}
    for any $s \leq t$ and any $\gamma \in L^2\big([0, t], \R^{d'\times q}\big)$.
\end{enumerate}
\end{lemma}

\begin{proof}
\begin{enumerate}
    \item Since $\E[\lVert X(t)\rVert^2]$ and $C(t)$ are locally bounded by Remark \ref{r:pp} and Lemma \ref{lem: p-rec}, there is for every $t \in \R_+$ some $K_t \geq 0$ such that $\E[\lVert a^c_\oo(s)+A^c_{\oo, :}(s)X(s)\rVert^2]\leq K_t$ and $\lVert C_\oo(s)\rVert \leq K_t$ for $0 \leq s \leq t$.
    Consequently, $\smash{\int_0^t\gamma(s)\dd X_\oo(s)}$ is well-defined for $\gamma\in L^2([0,t],\R^{d'\times q})$ and $\smash{\lVert \int_0^t\gamma(s)\dd X_\oo(s) \rVert_2 \leq K_t\lVert \gamma\rVert_{L^2([0,t])}}$. Now we can find simple functions $\smash{\gamma^{(n)} = \sum_{k=1}^{K_n} \gamma_k^{(n)} \bm{1}_{(t^{(n)}_{k-1}, t^{(n)}_{k}]}}$ such that $\lVert \gamma^{(n)} - \gamma \rVert_{L^2([0,t])} \to 0$. Thus $$\Big\lVert \int_0^t \gamma^{(n)}(s) \dd X_\oo(s) - \int_0^t \gamma(s) \dd X_\oo(s) \Big\rVert_2 \to 0.$$ Since $\int_0^t \gamma^{(n)}(s) \dd X_\oo(s) \in \mathcal{E}(X_\oo, t, d')$, this proves ``$\supseteq$'' in \eqref{eq: char_L2}. \\

    For the converse inclusion it suffices to show that the right-hand side of \eqref{eq: char_L2} is closed. To this end, let $F(t,s)$ denote the unique solution to $\frac{\partial}{\partial t}F(t,s) = A^c(t)F(t,s)$ with $F(s,s)=I_d$, which is locally bounded because $A^c$ is locally bounded. Applying \eqref{eq: X_closedform}, we get
    \begin{equation}\label{eq: Volterra_decomp}
        \int_0^t\gamma(r)\dd X_\oo(r) = U_\gamma + \int_0^t\Lambda_\gamma(s)\dd M(s),
    \end{equation}
    where $U_\gamma$ is a square-integrable, $\F_0$-measurable random variable and $\Lambda_\gamma(s) := \gamma(s)H + \int_s^t \gamma(r)A^c_{\oo, :}(r)F(r,s)\dd r$, with $H$ as defined at the beginning of this section. Define $$\xi_\gamma(s) := \Lambda_\gamma(s) C_{:, \oo}(s) C_\oo(s)^{-1} = \gamma(s) + \int_s^t\gamma(r)Q(r,s)\dd r,$$ where $Q(r,s) := A^c_{\oo, :}(r)F(r,s) C_{:, \oo}(s) C_\oo(s)^{-1}$, which is locally bounded by assumption. Hence the Volterra operator $(\mathcal V\gamma)(s) := \int_s^t\gamma(r)Q(r,s)\dd r$ is bounded on $L^2([0,t],\R^{d'\times q})$. Moreover, if $\lVert Q(r, s)\rVert \leq K_t$ for some $K_t \geq 0$ and $0 \leq s \leq r \leq t$, then iterated integration and Young's convolution inequality yield $\lVert \mathcal{V}^n\rVert \leq \frac{(K_t t)^n}{n!}$ and hence $\sum_{n=0}^\infty (-\mathcal V)^n$ converges in operator norm. In particular, $I+\mathcal V$ has a bounded inverse on $L^2([0,t])$. Therefore
    \begin{equation}\label{eq: Volterra_coercive}
        \|\gamma\|_{L^2([0,t])}
        \le K_t\|\xi_\gamma\|_{L^2([0,t])}
    \end{equation}
    for some $K_t \geq 0$. Moreover, $C_{:, \oo}(t) C_\oo^{-1}(t) C_{\oo, :}(t) \leq C(t)$, thus
    \begin{align*}
        \|\xi_\gamma(s)\|_F^2 &= \mathrm{tr}\left(\Lambda_\gamma(s)C_{:, \oo}(s) C_\oo^{-2}(s)C_{\oo, :}(s)\Lambda_\gamma(s)^\top\right)\\
        &\leq \|C_\oo(s)^{-1}\|\mathrm{tr}\left(\Lambda_\gamma(s)C_{:, \oo}(s) C_\oo^{-1}(s)C_{\oo, :}(s)\Lambda_\gamma(s)^\top\right)\\
        &\leq\|C_\oo(s)^{-1}\|\mathrm{tr}\left(\Lambda_\gamma(s)C(s)\Lambda_\gamma(s)^\top\right).
    \end{align*}
    Integrating both sides above, combining this with \eqref{eq: Volterra_coercive} and using the It\={o} isometry gives
    \begin{equation}\label{eq: gamma_coercive}
        \|\gamma\|_{L^2([0,t])}\leq K'_t \left\|\int_0^t\Lambda_\gamma(s)\dd M(s)\right\|_2
    \end{equation}
    for some $K'_t \geq 0$. If $V = \alpha+\beta X_\oo(0) +\int_0^t\gamma(s)\dd X_\oo(s),$ then $V = V_0+\int_0^t\Lambda_\gamma(s)\dd M(s)$ by \eqref{eq: Volterra_decomp} for some square-integrable, $\F_0$-measurable random variable $V_0$ and in particular $\big\|\int_0^t\Lambda_\gamma(s)\dd M(s)\big\|_2 \leq \|V\|_2$. Hence
    \begin{equation}\label{eq: final_coercive}
        \|\gamma\|_{L^2([0,t])} \leq K_t\|V\|_2 .
    \end{equation}
    Now let $V_n = \alpha_n+\beta_nX_\oo(0) +\int_0^t\gamma_n(s)\dd X_\oo(s)$ form a Cauchy sequence in $L^2(X, t, d')$. By \eqref{eq: final_coercive}, $(\gamma_n)_{n \in \N}$ is Cauchy in $L^2([0,t])$, hence $\gamma_n\to\gamma$ for some $\gamma \in L^2([0, t])$. In particular, $\int_0^t\gamma_n(s)\dd X_\oo(s) \to \int_0^t\gamma(s)\dd X_\oo(s)$ 
    in $L^2(X, t, d')$ because $\lVert \int_0^t\gamma(s)\dd X_\oo(s) \rVert_2 \leq K_t\lVert \gamma\rVert_{L^2([0,t])}$. Therefore $\widetilde L^2_0(X, t, d'):=\{\int_0^t \gamma(s) \dd X(s): \gamma \in L^2([0, t])\}$ is closed. A fortiori, the right-hand side of \eqref{eq: char_L2} is also closed as the sum of a finite-dimensional subspace and the closed subspace $\smash{\widetilde L^2_0(X, t, d')}$, see for example \textcite[Corollary 7-4.9]{Ma2002}.
    \item Define the linear map $J_t : \mathcal{E}(X, t, d') \to L^2(Y, t, d')$ by $\alpha + \sum_{i=1}^n \gamma_i X(t_i) \mapsto \alpha + \sum_{i=1}^n \gamma_i Y(t_i)$. For $Z \in \mathcal{E}(X, t, d')$ we have $\|J_t Z\| = \|Z\|$ because the first and second moments of $X$ and $Y$ match. Since $\mathcal{E}(X, t, d')$ is dense in $L^2(X, t, d')$, there exists a unique linear isometry $I_t : L^2(X, t, d') \to L^2(Y, t, d')$ such that $I_t$ and $J_t$ coincide on $\mathcal{E}(X, t, d')$. It remains to show that $I_t$ is surjective. Suppose $Z \in L^2(Y, t, d')$ so that $\lVert Z_n - Z \rVert_2 \to 0$ for some $Z_n \in \mathcal{E}(Y, t, d')$. By isometry we conclude that $H_n := I^{-1}_t(Z_n)$ is a Cauchy sequence converging to some $H \in L^2(X, t, d')$. Finally $\|I_t H - Z\| = \lim_{n \to \infty} \|Z_n - Z\| = 0$ and therefore $I_t(H) = Z$. More generally, for every $0\leq s\leq t$, we have $I_t\bigl(\mathcal E(X_\oo,s,d')\bigr)=\mathcal E(Y_\oo,s,d')$. Taking closures yields $I_t\bigl(L^2(X_\oo,s,d')\bigr) = L^2(Y_\oo,s,d')$. \\

    Finally, \eqref{eq: isometry} is evident if $\gamma$ is a simple function. For general $\gamma \in L^2\big([0, t], \R^{d'\times q}\big)$, choose a sequence of simple $\gamma^{(n)}$ with $\|\gamma^{(n)} - \gamma\|_{L^2([0, t])} \to 0$. By the proof of part 1, we have
    \[\Big\Vert \int_0^s \gamma^{(n)}(r) \dd X_\oo(r) - \int_0^s \gamma(r) \dd X_\oo(r) \Big\rVert_2 \to 0, \quad \Big\Vert \int_0^s \gamma^{(n)}(r) \dd Y_\oo(r) - \int_0^s \gamma(r) \dd Y_\oo(r) \Big\rVert_2 \to 0\]
    as $n \to \infty$. The claim then follows because the continuity of the operator $I_t$ yields
    \[I_t \int_0^s \gamma(r) \dd X_\oo(r) = \lim_{n \to \infty} I_t \int_0^s \gamma^{(n)}(r) \dd X_\oo(r) = \lim_{n \to \infty} \int_0^s \gamma^{(n)}(r) \dd Y_\oo(r) = \int_0^s \gamma(r) \dd Y_\oo(r). \]
    \end{enumerate}
\end{proof}

The following main results of this subsection now establish that the K\'alm\'an--Bucy filtering, prediction and smoothing equations from Propositions \ref{prop: KalBuc}\textendash \ref{prop: KalBuc_smooth} yield the optimal linear filter for polynomial processes of order 2:

\begin{corollary}[Filtering]\label{coro: Filt}
Suppose that $X$ is a polynomial process of order 2 as in \eqref{eq:pp1} and let $C(t) = \E(N^c(t))$ with $N^c(t)$ as in Lemma \ref{lem: p-rec}. Assume that $\mathrm{col} \:A^c_{\oo, \uu}(t) \subseteq \mathrm{col} \:C_\oo(t)$ holds for almost every $t \in \R_+$ and that $\smash{\int_0^t \lVert (C_\oo(s)^+)^{1/2} a^c_{\oo}(s)\rVert^2 \dd s< \infty}$ and $\smash{\int_0^t \lVert (C_\oo(s)^+)^{1/2} A^c_{\oo, :}(s) \dd s\rVert^2 < \infty}$ for every $t \in \R_+$. Let $\smash{\widehat \Sigma}$ and $\smash{\widehat{X}}$ denote the unique solutions of \eqref{eq: Kalman_eq_P} and \eqref{eq: Kalman_eq}. Then $\smash{\widehat{X}(t)}$ is an optimal linear filter for $X(t)$ and $$\E\bigl( \big(X(t) - \widehat X(t)\big) \smash{\big( X(t) - \widehat X(t)\big)}^\top\bigr) = \widehat{\Sigma}(t).$$
\label{Xenun34-4}
\end{corollary}
\begin{proof}
  Let $Y$ denote the Gaussian equivalent of $X$. By Lemma \ref{lem: funcana}.2, there exists an isometric isomorphism $I_t:L^2(X,t,d)\to L^2(Y,t,d)$ such that $I_t(X(t))=Y(t)$ and $I_t\bigl(L^2(X_\oo,t,d)\bigr) = L^2(Y_\oo,t,d).$ Consequently, $\|X(t)-Z\|_2 = \|Y(t)-I_t(Z)\|_2$ for $Z\in L^2(X_\oo,t,d)$. The linear filtering problem of minimising  $\mathbb{E}(\lVert X(t) - Z\rVert^2)$ over $Z \in L^2(X_\oo, t, d)$ is thus equivalent to minimising $\mathbb{E}(\lVert Y(t) - Z\rVert^2)$ over $Z \in L^2(Y_\oo, t, d)$, for which explicit solutions are given by the K\'alm\'an--Bucy filter from Proposition \ref{prop: KalBuc}. In particular, the optimal linear filter for $X$ is equal to the inverse image under $I_t$ of the optimal linear filter for $Y$ and hence given by the same K\'alm\'an--Bucy SDE from Proposition \ref{prop: KalBuc}. Indeed, writing $\widehat X^X$ and $\widehat X^Y$ for the solutions of \eqref{eq: Kalman_eq} driven by $X_\oo$ and $Y_\oo$, respectively,  we get $\widehat X^X(t) = g(X_\oo(0)) + \int_0^t h_t(s)\dd X_\oo(s)$ and $\widehat X^Y(t) = g(Y_\oo(0)) +  \int_0^t h_t(s)\dd Y_\oo(s)$ with some deterministic function $h_t$ and some affine function $g$. Approximating $h_t$ by deterministic step functions and using the definition and continuity of $I_t$ yields $I_t(\int_0^t h_t(s)\dd X_\oo(s))=\int_0^t h_t(s)\dd Y_\oo(s)$ and therefore $I_t(\widehat X^X(t))=\widehat X^Y(t)$. Hence the inverse image of the K\'alm\'an--Bucy filter of $Y$ is precisely the \mbox{solution $\widehat X(t)$ of \eqref{eq: Kalman_eq} driven by $X_\oo$.}
\end{proof}

\begin{corollary}[Prediction]
Suppose that $X$ is a polynomial process of order 2 as in Corollary \ref{coro: Filt}. For $s < t$ let $\widehat \Sigma(t, s)$ and $\widehat{X}(t, s)$ denote the unique solutions of \eqref{eq: Kalman_pred_P} and \eqref{eq: Kalman_pred}. Then $\widehat{X}(t, s)$ is an optimal linear predictor for $X(t)$ and $\E\bigl( \big(X(t) - \widehat X(t, s)\big) \smash{\big( X(t) - \widehat X(t, s)\big)^\top\bigr) = \widehat{\Sigma}(t, s)}$.
\label{Xenun35-5}
\end{corollary}
\begin{proof}
  This follows by the same arguments as in the proof of Corollary \ref{coro: Filt}.
\end{proof}

\begin{corollary}[Smoothing]
Suppose that $X$ is a polynomial process of order 2 as in Corollary \ref{coro: Filt} and assume
moreover that $\widehat\Sigma^+(t)$ from Corollary \ref{coro: Filt} is locally bounded. For $s > t$ let $\widehat \Sigma(t, s)$ and $\widehat{X}(t, s)$ be defined by Eqs. \eqref{eq: Kalman_smooth_P} and \eqref{eq: Kalman_smooth}. Then $\widehat{X}(t, s)$ is an optimal linear smoother for $X(t)$ and $\E\bigl( \big(X(t) - \widehat X(t, s)\big) \smash{\big( X(t) - \widehat X(t, s)\big)^\top\bigr) = \widehat{\Sigma}(t, s)}$.
\label{Xenun36-6}
\end{corollary}
\begin{proof}
  This follows by the same arguments as in the proof of Corollary \ref{coro: Filt}.
\end{proof}

\begin{remark}    As in the discrete case, the above results allow for simple error estimates for the optimal linear filter of a polynomial process. As before, the mean squared error (MSE) for the optimal linear filter is obtained via $\E(\lVert X(t) - \widehat X(t, s)\rVert^2) = \mathrm{tr}\: \widehat \Sigma(t, s)$.
\label{Xenun11-11}
\end{remark}

\section*{Acknowledgements}
The authors would like to thank two anonymous referees for their valuable comments.

\printbibliography

\end{document}